\documentclass[a4paper,11pt]{amsart}
\usepackage{amsmath,amsthm,amssymb}
\usepackage[mathscr]{eucal}
 \usepackage{cite}
\usepackage{upgreek}
\usepackage[bookmarks=false]{hyperref}
\usepackage{enumerate}
\usepackage{amsmath}
\usepackage{mathrsfs}
\usepackage{blindtext}
\usepackage{scrextend}
\usepackage{enumitem}
\addtokomafont{labelinglabel}{\sffamily}
\usepackage{color}
\usepackage[at]{easylist}

\makeatletter
\newcommand\footnoteref[1]{\protected@xdef\@thefnmark{\ref{#1}}\@footnotemark}
\makeatother

\numberwithin{equation}{section}

\theoremstyle{plain}
\newtheorem{main}{Theorem}
\newtheorem{mcor}[main]{Corollary}

\newtheorem{theorem}{Theorem}[section]

\newtheorem{lemma}[theorem]{Lemma}

\theoremstyle{definition}
\newtheorem{definition}[theorem]{Definition}
\newtheorem*{definition*}{Definition}

\newtheorem{remark}[theorem]{Remark}

\begin{document}

\title[Stability for product groups and property $(\tau)$]
{ Stability for product groups and property $(\tau)$}

\author[A. Ioana]{Adrian Ioana}
\address{Department of Mathematics, University of California San Diego, 9500 Gilman Drive, La Jolla, CA 92093, USA}
\email{aioana@ucsd.edu}

\thanks{The author was supported in part by NSF Career Grant DMS \#1253402 and NSF FRG Grant \#1854074.}
\begin{abstract}  We study the notion of permutation stability (or P-stability) for countable groups.
Our main result provides a wide class of  non-amenable product groups which are not P-stable. This class includes the product group $\Sigma\times\Lambda$, whenever $\Sigma$ admits a non-abelian free quotient and $\Lambda$ admits an infinite cyclic quotient. In particular, we obtain that the groups $\mathbb F_m\times\mathbb Z^d$ and $\mathbb F_m\times\mathbb F_n$ are not P-stable, for any integers $m,n\geq 2$ and $d\geq 1$.
This implies that P-stability is not closed under the direct product construction, which answers a question of Becker, Lubotzky and Thom. The proof of our main result relies on a construction of asymptotic homomorphisms from $\Sigma\times\Lambda$ to finite symmetric groups starting from sequences of finite index subgroups in $\Sigma$ and  $\Lambda$ with and without property $(\tau)$. Our method is sufficiently robust to show that the groups covered are not even flexibly P-stable, thus giving the first such non-amenable residually finite examples. 

\end{abstract}

\maketitle

\section{Introduction and statement of main results}

The notion of permutation stability has been developed in a series of works \cite{GR09,AP14,BLT18}.
A countable group $\Gamma$ is {\it stable in permutations} (or P-{\it stable}) if any ``almost homomorphism" from $\Gamma$ to a finite symmetric group is ``close" to a homomorphism. To make this precise, we endow the symmetric group Sym$(X)$ of any finite set $X$ with the normalized Hamming metric:
$$\text{d}_{\text{H}}(\sigma,\tau)=\frac{1}{|X|} |\{x\in X\mid \sigma(x)\not= \tau(x)\}|.$$

Hereafter, we will use the same formula to define the normalized Hamming distance between any maps $\sigma$ and $\tau$ with domain  (but not necessarily co-domain) equal to $X$. 
 \begin{definition}\label{P} 
 A sequence of maps $\sigma_n:\Gamma\rightarrow\text{Sym}(X_n)$, for some finite sets $X_n$, is called an {\it asymptotic homomorphism} 
 if  $\lim\limits_{n\rightarrow\infty}\text{d}_{\text{H}}(\sigma_n(gh),\sigma_n(g)\sigma_n(h))=0$, for every $g,h\in\Gamma$. 
The group $\Gamma$ is called \text{P}-{\it stable}\footnote{Definition \ref{P} agrees with the definitions of P-stability given in \cite{AP14} when $\Gamma$ is finitely presented and in \cite{BLT18} when $\Gamma$ is finitely generated, see Lemma \ref{equivalence}.} if for any asymptotic homomorphism $\sigma_n:\Gamma\rightarrow\text{Sym}(X_n)$, there exists a sequence of homomorphisms $\tau_n:\Gamma\rightarrow \text{Sym}(X_n)$ such that $\lim\limits_{n\rightarrow\infty}\text{d}_{\text{H}}(\sigma_n(g),\tau_n(g))=0$, for every $g\in\Gamma$.
\end{definition} 
More generally, one can define stability with respect to any class $\mathcal C$ of metric groups endowed with bi-invariant metrics (see \cite{AP14,AP17,CGLT17,Th17}). While this notion has only been formalized recently, in the case when $\Gamma=\mathbb Z^2$ and $\mathcal C$ consists of groups of matrices, the stability problem has been studied extensively in the literature.  Indeed, this problem is equivalent to the well-known question (posed in \cite{Ro69} for the normalized Hilbert-Schmidt norm and in \cite{Ha76} for the operator norm) of whether ``almost commuting" matrices are ``close" to commuting matrices.  The answer depends both on the groups of matrices considered and the norms chosen (see the introduction of \cite{AP14}). For instance, if $\mathcal C$ is the class of unitary groups $\{\text{U}(n)\mid n\in\mathbb N\}$,  then the answer is positive if one uses the normalized Hilbert-Schmidt norm \cite{HL08,Gl10} and negative if one uses the operator norm \cite{Vo83}. Recently, the stability problem with respect to unitary groups has been investigated for general countable groups $\Gamma$ in \cite{HS17,ESS18}  and for other matrix norms in \cite{CGLT17,LO18}.

At the same time, there has been a surge of 
interest in the study of permutation stability.
This started with the works of Glebsky and Rivera \cite{GR09} who observed that finite groups are P-stable\footnote{\label{note1}The results referenced here are stated in \cite{GR09} using the notion of stability in permutations for presentations of groups, see Remark \ref{equations}. In the form presented here, they follow from \cite{AP14}, where it was shown that stability is a group property, i.e., it is independent of the choice of the presentation.}, and of Arzhantseva and P\u{a}unescu \cite{AP14} who proved that abelian groups are P-stable (see \cite{BM18} for a quantitative approach to these results).
In \cite{BLT18}, Becker, Lubotzky and Thom obtained a characterization of P-stability for amenable groups in terms of invariant random subgroups, which they used to show that polycyclic groups and the Baumslag-Solitar  groups BS$(1,n)$  are P-stable. On the other hand, Becker and Lubotzky  \cite{BL18} proved that property ($\tau$) groups $\Gamma$ are not P-stable by removing one point from a set on which 
$\Gamma$ acts and deforming the action to get an almost action.
This motivated them to define two flexible variants of stability  (see Definition \ref{FP}). Subsequently, Lazarovich, Levit and Minsky proved that surface groups are flexibly P-stable \cite{LLM19}.

The study of P-stability is motivated in part by the longstanding problem of whether any countable group is sofic.  By an observation in \cite{GR09}, in order to find a non-sofic group, it is enough to find a group that is both P-stable and non-residually finite\footnoteref{note1}.
 We note that this point of view was used by De Chiffre, Glebsky, Lubotzky and Thom in their breakthrough work \cite{CGLT17}  to construct non-Frobenius-approximable groups. Very recently, Burton and Bowen proved that the existence of non-sofic groups would also follow from the flexible P-stability of PSL$_d(\mathbb Z)$  for $d\geq 5$  \cite{BB19}. 

The above results have led to a much better understanding of permutation stability, by providing several classes of P-stable and non-P-stable groups, as well as potential applications of this notion.
However, in spite of the progress made, the following basic question posed in \cite{BLT18} is open: 
is P-stability closed under direct products? 
While P-stability is clearly closed under free products, but not under the amalgamated free product or semi-direct product constructions by results in \cite{BL18}, 
 the situation remained unclear for direct products.
 
We settle this question in the negative here, by giving the first examples of P-stable groups whose direct product is not P-stable (see Corollary \ref{B} and the paragraph following it). 
This will be deduced as a consequence of our main result (Theorem \ref{A}) which provides a general criterion for non-P-stability of direct products of groups.
Moreover, our method of proof is sufficiently robust to address the flexible versions of P-stability introduced in \cite{BL18}, allowing to prove the following:

\begin{main}\label{A}
Let $\Sigma$ and $\Lambda$ be finitely generated groups. Assume that  $\Sigma$ admits a free non-abelian quotient and $\Lambda$ does not have property $(\tau)$. Then $\Sigma\times\Lambda$ is not very flexibly P-stable. \end{main}

Before presenting some examples of groups covered by Theorem \ref{A}, let us discuss the notions used in its statement and an equivalent formulation of it.

A countable group $\Lambda$ has {\it property $(\tau)$} if the quasi-regular representation of $\Lambda$ on $\bigoplus_{[\Lambda:\Delta]<\infty}\ell^2_0(\Lambda/\Delta)$, where $\Delta$ runs through all finite index subgroups of $\Lambda$ and $\ell^2_0(\Lambda/\Delta):=\ell^2(\Lambda/\Delta)\ominus\mathbb C{\bf 1}_{\Lambda/\Delta}$, does not have almost invariant vectors \cite{Lu94}. Property ($\tau)$ is a weaker version of property (T) which is satisfied by any irreducible lattice in a product of second countable, locally compact non-compact groups, at least one of which has property (T) \cite{LZ89}. In the opposite direction, any group admitting an infinite, residually finite amenable quotient group does not have property $(\tau)$ \cite{LW93,LZ03}.

\begin{definition}\label{FP}
A countable group $\Gamma$ is called 
{\it flexibly} \text{P}-{\it stable} if for any asymptotic homomorphism $\sigma_n:\Gamma\rightarrow\text{Sym}(X_n)$, there exist a sequence of finite sets $Y_n$ and homomorphisms $\tau_n:\Gamma\rightarrow \text{Sym}(Y_n)$ such that $X_n\subset Y_n$, for every $n$,  $\lim\limits_{n\rightarrow\infty}\text{d}_{\text{H}}(\sigma_n(g),\tau_n(g)_{|X_n})=0$, for every $g\in\Gamma$, and $\lim\limits_{n\rightarrow\infty}\frac{|Y_n|}{|X_n|}=1$.

The group $\Gamma$ is called {\it very flexibly} \text{P}-{\it stable} if for any asymptotic homomorphism $\sigma_n:\Gamma\rightarrow\text{Sym}(X_n)$, there exist a sequence of finite sets $Y_n$ and homomorphisms $\tau_n:\Gamma\rightarrow \text{Sym}(Y_n)$ such that $X_n\subset Y_n$, for every $n$, and $\lim\limits_{n\rightarrow\infty}\text{d}_{\text{H}}(\sigma_n(g),\tau_n(g)_{|X_n})=0$, for every $g\in\Gamma$.
\end{definition}
\begin{remark} A group $\Gamma$ is very flexibly P-stable if any asymptotic homomorphism is essentially obtained by restricting homomorphisms $\tau_n:\Gamma\rightarrow\text{Sym}(Y_n)$ to ``almost invariant" sets $X_n\subset Y_n$, i.e., such that $|\tau_n(g)X_n\triangle X_n|/|X_n|\rightarrow 0$, for every $g\in\Gamma$. If the sets $X_n$ are obtained by removing $o(|Y_n|)$ points from $Y_n$ (in which case they are trivially almost invariant), then  $\Gamma$ is flexibly P-stable.

It is unclear how much weaker these notions are than the (strict) notion of P-stability. On the one hand, P-stability coincides with flexible P-stability for amenable groups and flexible P-stability coincides with very flexible P-stability for property $(\tau)$ groups (see Lemma \ref{vs}). On the other hand, it is open whether property ($\tau$) groups can be flexibly P-stable and whether surface groups are  P-stable (see \cite{BL18,LLM19}). Moreover, while very flexible P-stability is inherited by subgroups of finite index (see Lemma \ref{finindex}), we do not know if this holds for P-stability or flexible P-stabillity. 
\end{remark}
Since very flexible P-stability passes to finite index subgroups, Theorem \ref{A} implies the following seemingly stronger statement: the product  between a large group and a group without property $(\tau)$ (and any group containing such a product as a finite index subgroup) is not very flexibly P-stable.
Recall that a group is called {\it large} if one of its finite index subgroups admits a non-abelian free quotient. 
By \cite{BP78} any finitely presented group with at least two more generators than relators is large; for more recent examples of large groups, see \cite{La07} and the references therein.

Theorem \ref{A} thus provides a wide class of groups, including the product of any large group and any group having an infinite, residually finite amenable quotient, which are not very flexibly P-stable. As an immediate consequence, we derive the following concrete examples of non-P-stable groups:

\begin{mcor}\label{B} The following groups are not very flexibly P-stable:
\begin{enumerate}
\item $\mathbb F_m\times\mathbb Z^d$, for every integers $m\geq 2$ and $d\geq 1$. 
\item $\mathbb F_m\times\mathbb F_n$, for every integers $m,n\geq 2$.
\item the Baumslag-Solitar group $\emph{BS}(m,n)$, 
for every integers $m,n$ with $|m|=|n|\geq 2$.
\item the braid group  \emph{B}$_n$  and pure braid group \emph{PB}$_n$, for every integer $n\geq 3$.
\end{enumerate}
\end{mcor}
Since free groups are obviously stable and abelian groups are stable by \cite{AP14},  (1) and (2) imply that P-stability is not closed under direct products, thus answering Becker, Lubotzky and Thom's question \cite{BLT18} in the negative. Moreover, we deduce that a direct product of P-stable groups need not even be very flexibly P-stable.  However, since the groups we treat are not amenable, this leaves open the question of whether the product of two P-stable amenable groups is P-stable \cite{BLT18}.  

In \cite[Example 7.3]{AP14} it was shown that BS$(m,n)$ is P-stable if $m=n=\pm 1$ but not P-stable if $|m|\not=|n|$ and $|m|,|n|\geq 2$, while \cite[Theorem 1.2 (ii)]{BLT18} established that BS$(1,n)$ is stable for every $n\in\mathbb Z$.
Part (3) of Corollary \ref{B} completes the classification of P-stability of the Baumslag-Solitar groups BS$(m,n)$ by addressing the remaining case when $|m|=|n|\geq 2$.

To put Corollary \ref{B} into a better perspective, let us indicate several additional consequences of it.
First, as remarked in \cite[Section 4.4]{BL18} (extending observations made in \cite{GR09,AP14}), any group which is sofic and non-residually finite 
is not very flexibly P-stable. 
By \cite[Theorem 1.2 (iii)]{BLT18}, there are amenable residually finite groups which are not P-stable and thus not flexibly P-stable.
 Corollary \ref{B} gives the first examples of  non-amenable residually finite groups that are not flexibly P-stable, and of residually finite groups that are not very flexibly P-stable.



\begin{remark}\label{BB}
A countable group $\Gamma$ is called {\it Hilbert-Schmidt stable} (or HS-{\it stable})  if it is stable with respect to the class of unitary groups $\{(\text{U}(n),\text{d}_{\text{HS}})\mid n\in\mathbb N\}$ endowed with the normalized Hilbert-Schmidt distance $\text{d}_{\text{HS}}(T,S)=\|T-S\|_{\text{HS}}$ for $T,S\in\text{U}(n)$, where  $\|V\|_{\text{HS}}=\sqrt{n^{-1}\text{Tr}(V^*V)}$. 
 Since the normalized Hamming distance can be expressed in terms of the normalized Hilbert-Schmidt, the study of P-stability and HS-stability are similar in flavor \cite{AP14}. 

In spite of the similarity between these notions, Corollary \ref{B} highlights a surprising difference between them, by providing, to our knowledge, the first examples of HS-stable groups which are not P-stable.
By  \cite[Theorem 1]{HS17}, the product of two HS-stable groups is HS-stable provided that one of the groups is abelian (by \cite[Corollary D]{IS19} the same holds if one of the groups is amenable). Consequently,  $\mathbb F_m\times\mathbb Z^d$ is HS-stable but not P-stable, for any integers $m\geq 2$ and $d\geq 1$.


Note that is an open question whether HS-stability is closed under direct products. It seems likely that this question 
has a negative answer, and moreover that $\mathbb F_m\times\mathbb F_n$ is not HS-stable, for $m,n\geq 2$. Supporting evidence is provided by  \cite[Theorem E]{IS19} which shows that $\mathbb F_m\times\mathbb F_n$ is not stable with respect to the class $\{(\text{U}(M), \|\cdot\|_{2})\mid \text{$(M,\uptau)$ tracial von Neumann algebra} \}$ of unitary groups of tracial von Neumann algebras endowed with their $2$-norms, $\|T\|_2=\sqrt{\uptau(T^*T)}$.

\end{remark}

\begin{remark}\label{equations}
Let  $R\subset\mathbb F_k$ be a finite set, for $k\in\mathbb N$. The system of equations ($\star$) $r(\tau_1,...,\tau_k)=e$, for every $r\in R$, is called P-stable if  for every $\varepsilon>0$, there is $\delta>0$ such that the following holds: for any finite set $X$ and $\sigma_1,...,\sigma_k\in\text{Sym}(X)$ satisfying $\text{d}_{\text{H}}(r(\sigma_1,...,\sigma_k),\text{Id}_X)<\delta$, for every $r\in R$,  $(\star$) has a solution $\tau_1,...,\tau_k\in\text{Sym}(X)$ such that $\text{d}_{\text{H}}(\sigma_i,\tau_i)\leq\varepsilon$, for every $1\leq i\leq k$ (see \cite{GR09,AP14}). A finitely presented group $\Gamma=\langle\mathbb F_k | R\rangle$ is P-stable if and only if $R$ is P-stable \cite{AP14}. The P-stability of $\mathbb Z^2$ proved in \cite{AP14} thus implies  P-stability of the system $[a,b]=aba^{-1}b^{-1}=e$. 

On the other hand, since the groups $\mathbb F_2\times\mathbb Z$ and $\mathbb F_2\times\mathbb F_2$  are not P-stable by Corollary \ref{B}, we conclude that the systems $[a_1,b]=[a_2,b]=e$ and $[a_1,b_1]=[a_1,b_2]=[a_2,b_1]=[a_2,b_2]=e$ are not P-stable.
\end{remark}

Corollary \ref{B} also implies the existence of universal sofic groups which fail a certain lifting property for commuting subgroups. 
Let $\mathcal U$ be a free ultrafilter on $\mathbb N$ and $(X_n)$ finite sets with $\lim\limits_{n\rightarrow\mathcal U}|X_n|=+\infty$. Define the metric ultraproduct group $\prod_{\mathcal U}\text{Sym}(X_n):=\big(\prod_{n}\text{Sym}(X_n)\big)/\mathcal N$, where $\mathcal N$ is the subgroup of $(\sigma_n)\in\prod_n\text{Sym}(X_n)$ satisfying $\lim\limits_{n\rightarrow\mathcal U}\text{d}_{\text{H}}(\sigma_n,\text{Id}_{X_n})=0$. Since a countable group is sofic if and only if it embeds into  $\prod_{\mathcal U}\text{Sym}(X_n)$ \cite{ES04}, the latter  is called a {\it universal sofic group}.
\begin{mcor}
There exist  countable commuting subgroups $\Sigma, \Lambda$ of a universal sofic group $\prod_{\mathcal U}\emph{Sym}(X_n)$ such that the following holds: there are no commuting subgroups $\Sigma_n, \Lambda_n$ of $\emph{Sym}(X_n$), for all $n\in\mathbb N$, such that $\Sigma\subset\prod_{\mathcal U}\Sigma_n$ and $\Lambda\subset\prod_{\mathcal U}\Lambda_n$. 

\end{mcor}
We end the introduction by discussing a weakening of the notion of P-stability found by considering asymptotic homomorphisms that are sofic approximations \cite{AP14}. Let $\Gamma$ be a countable group.

\begin{definition}
An asymptotic homomorphism $\sigma_n:\Gamma\rightarrow\text{Sym}(X_n)$ is called a {\it sofic approximation} of $\Gamma$ if $\lim\limits_{n\rightarrow\infty}\text{d}_{\text{H}}(\sigma_n(g),\text{Id}_{X_n})=1$, for every $g\in\Gamma\setminus\{e\}$. The group $\Gamma$ is called {\it weakly \emph{P}-stable} (respectively, {\it weakly flexibly \emph{P}-stable} or {\it weakly very flexibly \emph{P}-stable}) if  the condition from Definition \ref{P} (respectively, the conditions from Definition \ref{FP}) holds for any sofic approximation $(\sigma_n)$ of $\Gamma$.
\end{definition}

The notion of weak P-stability is strictly weaker than that of P-stability.  More precisely, \cite[Theorem 7.2]{AP14} shows that any finitely presented, residually finite amenable group is weakly P-stable, whereas \cite[Theorem 1.2 (iii)]{BLT18} proves that there is such a group which is not P-stable.

Our last main result provides a class of non-amenable groups which are not weakly P-stable:
\begin{main}\label{C}
Any group which has a subgroup of finite index isomorphic to $\mathbb F_m\times\mathbb Z^d$ or to $\mathbb F_m\times\mathbb F_n$, for some integers $m,n\geq 2$ and $d\geq 1$, is not weakly very flexibly \emph{P}-stable. 
In particular,  any group from Corollary \emph{\ref{B}}, parts (1)-(3),  is not weakly very flexibly \emph{P}-stable.
\end{main}

Consequently, the Baumslag-Solitar BS$(m,n)$  group is not weakly P-stable,  whenever $|m|=|n|\geq 2$. This settles a question posed by Arzhantseva and P\u{a}unescu in \cite[Example 7.3]{AP14}. 
As a special case of Theorem \ref{C}, we deduce that $\mathbb F_2\times\mathbb Z$ is not weakly flexibly P-stable. This answers a question raised by Bowen and Burton in \cite{BB19} who emphasized that $\mathbb F_2\times\mathbb Z$ seems to be the most elementary group for which weak flexible P-stability was unknown (note that the notion of flexible stability used in \cite{BB19} is what we call here weak flexible stability). 

\subsection*{Comments on the proof of Theorem \ref{A}} We end the introduction with an outline of the proof of Theorem \ref{A} under the following additional assumption: there exist a group $\Gamma$, a sequence $\{\Gamma_n\}_{n=1}^{\infty}$ of finite index normal subgroups of $\Gamma$, and homomorphisms $q_n:\Lambda\rightarrow \Gamma/\Gamma_n$ such that 
\begin{itemize}
\item $\Sigma=\Gamma*\mathbb Z$, 
\item $\Gamma$ has property $(\tau)$ with respect to $\{\Gamma_n\}_{n=1}^{\infty}$, and
\item $\Lambda$ does not have property $(\tau)$ with respect to $\{\ker(q_n)\}_{n=1}^{\infty}$.
\end{itemize}
 
This assumption holds for $\Sigma=\mathbb F_3$ and $\Lambda=\mathbb Z$, by taking $\{\Gamma_n\}_{n=1}^{\infty}$ be a sequence  of finite index normal subgroups of $\Gamma=\mathbb F_2$ with property $(\tau)$ and $q_n:\Lambda\rightarrow\Gamma/\Gamma_n$ homomorphisms with $|q_n(\Lambda)|\rightarrow +\infty$. More generally, we use Kassabov's theorem \cite{Ka05} (that the symmetric groups $\{\text{Sym}(n)\}_{n=1}^{\infty}$ admit Cayley graphs which form a bounded degree expander family) to conclude that there is  $L\geq 2$ such that the assumption is satisfied when $\Sigma=\mathbb F_{L+1}$, $\Gamma=\mathbb F_L$ and $\Lambda$ is any group without property $(\tau)$.
 
Next,  let $X_n=\Gamma/\Gamma_n$, $p_n:\Gamma\twoheadrightarrow X_n$ be the quotient homomorphism. 
 and view $\Gamma\times\Lambda$ as a subgroup of $\Sigma\times\Lambda$. We define the left-right multiplication action $\sigma_n:\Gamma\times\Lambda\rightarrow\text{Sym}(X_n)$ by letting $$\text{$\sigma_n(g,h)x=p_n(g)xq_n(h)^{-1}$, for every $g\in\Gamma, h\in\Lambda, x\in X_n$. }$$
 
 There are two main ingredients in the proof of Theorem \ref{A}.
 
 The first is a rigidity result for asymptotic homomorphisms $\widetilde\sigma_n:\Sigma\times\Lambda\rightarrow\text{Sym}(X_n)$ extending $\sigma_n$, i.e., ${\widetilde{\sigma_n}}_{|\Gamma\times\Lambda}=\sigma_n$.
  Assume there are homomorphisms $\tau_n:\Sigma\times\Lambda\rightarrow\text{Sym}(Y_n)$, with $Y_n\supset X_n$ finite, such that $\text{d}_{\text{H}}(\widetilde\sigma_n(g),\tau_n(g)_{|X_n})\rightarrow 0$, for all $g\in\Sigma\times\Lambda$. Using the property $(\tau)$ assumption, we prove that there must be homomorphisms $\overline{\sigma}_n:\Sigma\times\Lambda\rightarrow\text{Sym}(X_n)$ extending $\sigma_n$ such that $\text{d}_{\text{H}}(\widetilde\sigma_n(g),\overline\sigma_n(g))\rightarrow 0$, for all $g\in\Sigma\times\Lambda$ (see Theorem \ref{conjugation}).
In other words, if $\widetilde\sigma_n$ is close to (the restriction to $X_n$ of) a homomorphism, then $\widetilde\sigma_n$ is close to a homomorphism which extends $\sigma_n$.
 
 The second ingredient in the proof of Theorem \ref{A} is the construction of a ``non-trivial" asymptotic homomorphism $\widetilde\sigma_n:\Sigma\times\Lambda\rightarrow\text{Sym}(X_n)$ extending $\sigma_n$. 
Using that $\Lambda$ does not have property $(\tau)$ with respect to $\{\ker(q_n)\}_{n=1}^{\infty}$, we construct in Lemma \ref{tech2} a permutation $\rho_n\in\text{Sym}(X_n)$ such that \begin{enumerate}
 \item $\text{d}_{\text{H}}(\rho_n\circ\sigma_n(e,h), \sigma_n(e,h)\circ\rho_n)\rightarrow 0$, for every $h\in\Lambda$, and
 \item $\max\{\text{d}_{\text{H}}(\rho_n\circ \sigma_n(e,h), \sigma_n(e,h)\circ\rho_n)\mid h\in\Lambda\}\geq \frac{1}{126}$, for infinitely many $n$.
 \end{enumerate}
Specifically, we first find $A_n\subset X_n$ which is almost invariant under the right multiplication action of $\Lambda$ and satisfies $\frac{|A_n|}{|X_n|}\in (\frac{1}{7},\frac{1}{6})$ for $n$ large (see Lemma \ref{AE}). After replacing $A_n$ with a subset, we may assume that $A_n\cap g_nA_n=\emptyset$, for $g_n\in X_n$. We then show that $\rho_n$ defined by $\rho_n(x)=g_nx$ if $x\in A_n$, $\rho_n(x)=g_n^{-1}x$ if $x\in g_nA_n$, and $\rho_n(x)=x$ if $x\notin A_n\cup g_nA_n$, satisfies conditions (1) and (2). 
 
Finally, condition (1) allows us to define an asymptotic homomorphism $\widetilde\sigma_n:\Sigma\times\Lambda\rightarrow\text{Sym}(X_n)$ which extends $\sigma_n$ by letting $\widetilde\sigma_n(t,e)=\rho_n$, where $t\in\mathbb Z$ is a generator. On the other hand, (2) guarantees that $\widetilde\sigma_n$ is not close to any homomorphism $\overline{\sigma}_n:\Sigma\times\Lambda\rightarrow\text{Sym}(X_n)$ which extends $\sigma_n$. But then the first ingredient above implies that $\Sigma\times\Lambda$ is not very flexibly P-stable, as desired.
 
In the general case, when $\Sigma$ is only assumed to have a non-abelian free quotient, after replacing it with a finite index subgroup, we may assume that there is an onto homomorphism $\pi:\Sigma\rightarrow\mathbb F_{L+1}$. Let $\widetilde\sigma_n:\mathbb F_{L+1}\times\Lambda\rightarrow\text{Sym}(X_n)$ be the asymptotic homomorphism constructed above which witnesses that $\mathbb F_{L+1}\times\Lambda$ is not very flexibly P-stable. Then we analyze the asymptotic homomorphism $\widetilde\sigma_n\circ (\pi\times\text{Id}_{\Lambda}):\Sigma\times\Lambda\rightarrow\text{Sym}(X_n)$ to show that $\Sigma\times\Lambda$ is not very flexibly P-stable.

\subsection*{Acknowledgements} I would like to thank Goulnara Arzhantseva and Pieter Spaas for several helpful comments and corrections, and Lewis Bowen and Andreas Thom for stimulating discussions.

\section{Preliminaries}
In this section, we first recall some notation and then gather several results that will be needed later. 
Let $X$ be a finite set. We denote by $\text{B}(\ell^2(X))$ the algebra of all linear maps $T:\ell^2(X)\rightarrow\ell^2(X)$ and by $\{\delta_x\}_{x\in X}$ the usual orthonormal basis of $\ell^2(X)$.

 The {\it normalized Hilbert-Schmidt norm} of $T\in\text{B}(\ell^2(X))$ is given by $$
\|T\|_{\text{HS}}=\sqrt{\frac{1}{|X|}\text{Tr}(T^*T)}=\sqrt{\frac{1}{|X|}\sum_{x,y\in X}|\langle T\delta_x,\delta_y\rangle|^2}.$$

Let $U:\text{Sym}(X)\rightarrow\text{U}(\ell^2(X))$ be the group homomorphism given by $U_\sigma(\delta_x)=\delta_{\sigma(x)}$, for all $x\in X$.
Hereafter, we view $\text{Sym}(X)$ as a subgroup of U$(\ell^2(X))$, via the embedding $U$.
Note that $$\text{$\|U_{\sigma}-U_{\tau}\|_{\text{HS}}=\sqrt{2\;\text{d}_{\text{H}}(\sigma,\tau)}$, for every $\sigma,\tau\in\text{Sym}(X)$.}$$

\subsection{On the distance to invariant sets}  Next, we record the following well-known fact.
\begin{lemma}\label{component}
Let $Y$ be a finite set, $X\subset Y$ be a subset and $H<\emph{Sym}(Y)$ be a subgroup. Then  there exists an $H$-invariant subset $X_0\subset Y$ such that $|X_0\triangle X|\leq 2\max_{h\in H}|X\triangle hX|$.
\end{lemma}

{\it Proof.} Put $\varepsilon=\max_{h\in H}|X\triangle hX|$ and define the $H$-invariant function $f=\frac{1}{|H|}\sum_{h\in H}{\bf 1}_{hX}\in\ell^1(Y)$.
Since $\|{\bf 1}_X-{\bf 1}_{hX}\|_1=|X\triangle hX|\leq\varepsilon$, for every $h\in H$, we get that $\|{\bf 1}_X-f\|_{1}\leq \varepsilon$. Then the set $X_0=\{y\in Y\mid f(y)\geq \frac{1}{2}\}$ is $H$-invariant and since $$\|{\bf 1}_X-f\|_1=\sum_{y\in Y\setminus X}|f(y)|+\sum_{y\in X}|f(y)-1|\geq \frac{1}{2}|(Y\setminus X)\cap X_0|+\frac{1}{2}|X\setminus X_0|=\frac{1}{2}|X_0\triangle X|,$$
the conclusion follows. \hfill$\blacksquare$

\subsection{Kazhdan constants}
We continue by recalling the notion of a Kazhdan constant and two well-known facts which we prove for completeness.
\begin{definition}
Let $G$ be a finite group and $S$ be a set of generators. The {\it  Kazhdan constant} $\kappa(G,S)$ is the largest constant $\kappa>0$ such that $\kappa\;\|\xi\|\leq\max_{g\in 
S}\|\pi(g)\xi-\xi\|$, for every $\xi\in \mathcal H$ and unitary representation $\pi:G\rightarrow \text{U}(\mathcal H)$ on a Hilbert space $\mathcal H$ without non-zero invariant vectors.
\end{definition}

\begin{lemma}\label{expansion}
Let $G$ be a finite group and $S$ be a set of generators. Then for every subset $A\subset G$ we have that $\kappa(G,S)^2\; |A| \;|G\setminus A|\leq\max_{g\in S}|gA\triangle A|\; |G|$.
\end{lemma}

{\it Proof.} Let $\lambda:G\rightarrow\text{U}(\ell^2(G))$ be the left regular representation. Put $\xi={\bf 1}_A-\frac{|A|}{|G|}{\bf 1}_G\in\ell^2(G)\ominus\mathbb  C{\bf 1}_G$. Then the conclusion  is equivalent to the inequality $\kappa(G,S)\;\|\xi\|_2\leq\max_{g\in S}\|\lambda(g)\xi-\xi\|_2$, which holds since the restriction of $\lambda$ to $\ell^2(G)\ominus\mathbb C{\bf 1}_G$ has no non-zero invariant vectors. \hfill$\blacksquare$

\begin{lemma}\label{almostinv}
Let $G$ be a finite group and $S$ be a set of generators. Then for every unitary representation $\pi:G\rightarrow\emph{U}(\mathcal H)$ and $\xi\in\mathcal H$
we have that $\kappa(G,S)\cdot\max_{g\in G}\|\pi(g)\xi-\xi\|\leq 2\;\max_{g\in S}\|\pi(g)\xi-\xi\|$.
\end{lemma}

{\it Proof.} Let $\mathcal H^{G}$ be the subspace of $\mathcal H$ consisting of $\pi(G)$-invariant vectors.
Let $\xi\in \mathcal H$ and write $\xi=\xi_1+\xi_2$, where $\xi_1\in \mathcal H\ominus\mathcal H^G$ and $\xi_2\in \mathcal H^G$. Then $\|\pi(g)\xi-\xi\|=\|\pi(g)\xi_1-\xi_1\|\leq 2\|\xi_1\|$, for every $g\in G$. Since the restriction of $\pi$ to $\mathcal H\ominus\mathcal H^{G}$ has no non-zero invariant vectors, we get that $\kappa(G,S)\;\|\xi_1\|\leq\max_{g\in S}\|\pi(g)\xi_1-\xi_1\|$ and the conclusion follows.
 \hfill$\blacksquare$

\subsection{Property $(\tau)$}
We are now ready to recall an equivalent formulation of property $(\tau$)  with respect to a sequence of finite index normal subgroups \cite{Lu94} .
\begin{definition}\label{tau}
Let $\Gamma$ be a finitely generated group and $S$ be a finite set of generators. Then $\Gamma$ has {\it property ($\tau$)} with respect of a sequence of finite index subgroups $\{\Gamma_n\}_{n=1}^{\infty}$ if $\inf_n\kappa(\Gamma/\Gamma_n,p_n(S))>0$, where $p_n:\Gamma\rightarrow\Gamma/\Gamma_n$ denotes the quotient homomorphism.
\end{definition}

If $\lim\limits_{n\rightarrow\infty}\kappa(\Gamma/\Gamma_n,p_n(S))=0$, then there exist sets $C_n\subset\Gamma/\Gamma_n$ satisfying $0<|C_n|<|\Gamma/\Gamma_n|/2$ which are almost invariant, in the sense that $\lim\limits_{n\rightarrow\infty}|p_n(g)C_n\triangle C_n|/|C_n|=0$, for every $g\in\Gamma$ (see \cite[Proposition 2.5]{LZ03}). 
Moreover,  Ab\'{e}rt and Elek  \cite[Theorem 4]{AE10} proved that one can choose $C_n$ such that the sequence $\{|C_n|/|\Gamma/\Gamma_n|\}_{n=1}^{\infty}$ converges to any prescribed limit in $[0,\frac{1}{2}]$. 

The next lemma, which is of independent interest and will be used in the proof of Lemma \ref{tech2}, generalizes this result to arbitrary, not necessarily decreasing, sequences of normal subgroups.  
\begin{lemma}
\label{AE}
In the notation of Definition \ref{tau}, assume that  $\lim\limits_{n\rightarrow\infty}\kappa(\Gamma/\Gamma_n,p_n(S))=0$ and $\Gamma_n<\Gamma$ is normal, for any $n$. Let $0<\alpha<\beta\leq\frac{1}{2}$.
Then for  large enough $n$ there is $C_n\subset\Gamma/\Gamma_n$ such that  $$\text{$\alpha\leq\frac{|C_n|}{|\Gamma/\Gamma_n|}\leq\beta$ and $\lim\limits_{n\rightarrow\infty}\frac{|p_n(g)C_n\triangle C_n|}{|\Gamma/\Gamma_n|}=0$, for every $g\in\Gamma$.}$$
\end{lemma}

{\it Proof.} If $\{\Gamma_n\}_{n=1}^{\infty}$ is a descending chain, the lemma is a direct consequence of \cite[Theorem 4]{AE10}.
In general, denote $G_n=\Gamma/\Gamma_n$ for $n\geq 1$. The proof is based on the following:

{\bf Claim}. Let $D_n\subset G_n$ be a sequence of sets such that 
$\lim\limits_{n\rightarrow\infty}\frac{|p_n(g)D_n\triangle D_n|}{|D_n|}=0$ and $0<|D_n|<\frac{3|G_n|}{4}$, for every $g\in\Gamma$ and $n\geq 1$.  Then for any large enough $n$ we can find $h_n\in G_n$ such that \begin{equation} \label{h_n}\frac{|D_n|^2}{4|G_n|}\leq |D_nh_n\cap D_n|\leq\frac{3|D_n|}{4}.\end{equation}

{\it Proof of the claim.}
Assume that the claim is false. After passing to a subsequence, we may assume that for every $n\geq 1$ and $h\in G_n$ we have $|D_nh\cap D_n|<\frac{|D_n|^2}{4|G_n|}$ or $|D_nh\cap D_n|>\frac{3|D_n|}{4}$. Let $H_n$ be the set of $h\in G_n$ such that $|D_nh\cap D_n|>\frac{3|D_n|}{4}$. If $h,h'\in H_n$, then $|D_nhh'\cap D_n|>\frac{|D_n|}{2}>\frac{|D_n|^2}{4|G_n|}$ and hence $hh'\in H_n$. This implies that $H_n$ is a subgroup of $G_n$. Next, since $$|D_n|^2=\sum_{h\in G_n}|D_nh\cap D_n|=\sum_{h\in H_n}|D_nh\cap D_n|+\sum_{h\in G_n\setminus H_n}|D_nh\cap D_n|\leq |D_n|\; |H_n|+\frac{|D_n|^2}{4|G_n|}\; |G_n|,$$ we get that $|H_n|\geq\frac{3|D_n|}{4}$. On the other hand, since $$\sum_{x\in D_n}|H_n\cap x^{-1}D_n|=\sum_{h\in H_n}|D_nh\cap D_n|\geq \frac{3|D_n|}{4}\; |H_n|,$$ we can find $x_n\in D_n$ such that $|x_nH_n\cap D_n|=|H_n\cap x_n^{-1}D_n|\geq\frac{3|H_n|}{4}.$ In particular, $|D_n|\geq\frac{3|H_n|}{4}$. 
Since $|D_n|\leq\frac{4|H_n|}{3}$, we get $|x_nH_n\triangle D_n|=|D_n|+|H_n|-2|x_nH_n\cap D_n|\leq |D_n|-\frac{|H_n|}{2}\leq\frac{5|H_n|}{6}.$ Thus, for every $g\in\Gamma$ we have
$$ |p_n(g)x_nH_n\triangle x_nH_n|\leq 2|x_nH_n\triangle D_n|+|p_n(g)D_n\triangle D_n|\leq\frac{5|H_n|}{3}+|p_n(g)D_n\triangle D_n|.$$
Since $\lim\limits_{n\rightarrow\infty}\frac{|p_n(g)D_n\triangle D_n|}{|D_n|}=0$ and $\frac{|D_n|}{|H_n|}\leq\frac{4}{3}$, it follows that $\limsup_{n\rightarrow\infty}\frac{|p_n(g)x_nH_n\triangle x_nH_n|}{|H_n|}\leq\frac{5}{3}<2$. Thus, for every $g\in\Gamma$ we have $p_n(g)\in x_nH_nx_n^{-1}$, for $n$ large enough. Since $\Gamma$ is finitely generated, we get that $H_n=G_n$, for $n$ large enough. This contradicts that $|H_n|\leq\frac{4|D_n|}{3}<|G_n|$, for any $n$.
\hfill$\square$

 Now, let $L$ be the set of $\ell\in [0,\frac{1}{2}]$ for which there is a sequence of nonempty sets $D_n\subset G_n$ with 
 $$\text{$\limsup_{n\rightarrow\infty}\frac{|D_n|}{|G_n|}=\ell$\;\; and \;\;$\lim\limits_{n\rightarrow\infty}\frac{|p_n(g)D_n\triangle D_n|}{|D_n|}=0$, for every $g\in\Gamma$.}$$ Since  $\lim\limits_{n\rightarrow\infty}\kappa(G_n,p_n(S))=0$ we have that $L\not=\emptyset$  (see, e.g., \cite[Proposition 2.5]{LZ03}). 
 
 We claim that $\inf L=0$. 
If $0\in L$, there is nothing to prove. Otherwise, let $\ell\in L\setminus\{0\}$ and $D_n\subset G_n$ be sets witnessing that $\ell\in L$. By the above claim for every $n$ large enough we can find $h_n\in G_n$ such that $\frac{|D_n|^2}{4|G_n|}\leq |D_nh_n\cap D_n|\leq\frac{3|D_n|}{4}.$
For every $n\geq 1$, define $D_n'=\begin{cases}\text{$D_nh_n\cap D_n$, if $\frac{|D_n|}{|G_n|}>\frac{\ell}{2}$}\\\text{$D_n$, if $\frac{|D_n|}{|G_n|}\leq\frac{\ell}{2}$.}\end{cases}$ If $\frac{|D_n|}{|G_n|}>\frac{\ell}{2}$, then  $p_n(g)D_n'\triangle D_n'\subset (p_n(g)D_n\triangle D_n)\cup (p_n(g)D_n\triangle D_n)h$ and hence we get that $$\frac{|p_n(g)D_n'\triangle D_n'|}{|D_n'|}\leq\frac{2|p_n(g)D_n\triangle D_n|}{\frac{|D_n|^2}{4|G_n|}}\leq \frac{16}{\ell}\;\frac{|p_n(g)D_n\triangle D_n|}{|D_n|}.$$
From this it follows that $\lim\limits_{n\rightarrow\infty}\frac{|p_n(g)D_n'\triangle D_n'|}{|D_n'|}=0$, for every $g\in\Gamma$. Thus, $\ell'=\limsup_{n\rightarrow\infty}\frac{|D_n'|}{|G_n|}\in L$. 
Since $\frac{|D_n'|}{|G_n|}\leq\max\{\frac{3|D_n|}{4|G_n|},\frac{\ell}{2}\}$, for every $n$, we conclude that $\ell'\leq\frac{3\ell}{4}$. This implies that $\inf L=0$.

Let now $0<\alpha<\beta\leq\frac{1}{2}$. Since $\inf L=0$, we can find a sequence of sets $D_n\subset G_n$ such that $\frac{|D_n|}{|G_n|}\leq\min\{\beta-\alpha,\alpha\}$, for $n$ large enough, and $\lim\limits_{n\rightarrow\infty}\frac{|p_n(g)D_n\triangle D_n|}{|D_n|}=0$, for every $g\in\Gamma$. 

For $n\geq 1$, let $k_n=\Big{\lceil}\frac{\log(1-\alpha)}{\log(1-\frac{|D_n|}{|G_n|})}\Big{\rceil}$ be the smallest integer such that $1-\big(1-\frac{|D_n|}{|G_n|}\big)^{k_n}\geq\alpha$. 
Let $m_n\geq 1$ be the smallest integer for which there exists a set $F_n\subset G_n$ of cardinality $m_n${ such that $C_n:=D_nF_n$ satisfies $\frac{|C_n|}{|G_n|}\geq\alpha$.
By \cite[Lemma 2.3]{AE10} we have that $m_n\leq k_n$. 
Then $\frac{|C_n|}{|G_n|}<\beta$, for all $n$. Indeed, 
if $g\in F_n$, then the minimality of $m_n$ implies that $\frac{|D_n(F_n\setminus\{g\})|}{|G_n|}<\alpha$ and thus $$\frac{|C_n|}{|G_n|}\leq\frac{|D_n(F_n\setminus\{g\})|}{|G_n|}+\frac{|D_ng|}{|G_n|}<\alpha+(\beta-\alpha)=\beta.$$

Finally, if $g\in\Gamma$, then $p_n(g)C_n\triangle C_n\subset\cup_{h\in F_n}(p_n(g)D_nh\triangle D_nh)$ and thus \begin{align*}\frac{|p_n(g)C_n\triangle C_n|}{|G_n|}\leq\frac{m_n\;|p_n(g)D_n\triangle D_n|}{|G_n|}\leq\frac{k_n|D_n|}{|G_n|}\;\frac{|p_n(g)D_n\triangle D_n|}{|D_n|}\end{align*} Since the sequence $\{\frac{k_n|D_n|}{|G_n|}\}_{n=1}^{\infty}$ is bounded, this implies that $\lim\limits_{n\rightarrow\infty}\frac{|p_n(g)C_n\triangle C_n|}{|C_n|}=0$, for every $g\in\Gamma$, which finishes the proof of the lemma.
\hfill$\blacksquare$

\section{Basic results on P-stability} In this section, we record three results on the general theory of P-stability. Note that with one exception, Lemma \ref{finindex}, these results will not be needed in the rest of the paper. 

\subsection{Equivalence of definitions of P-stability} The notion of P-stability was introduced in \cite[Definition 3.2]{AP14} (see also \cite{GR09}) for finitely presented groups, and generalized  to finitely generated groups in \cite[Definition 3.11]{BLT18}. Our next result provides an equivalent formulation of P-stability, in the sense of Definition \ref{P},  for general groups. This implies that for finitely generated groups the notions of P-stability given by \cite[Definition 3.11]{BLT18} and Definition \ref{P} coincide. 

Let $\Gamma$ be a countable group and $S$ a set of generators. Denote by $\{\bar{s}\}_{s\in S}$ the free generators of $\mathbb F_S$  and by $\pi:\mathbb F_S\rightarrow\Gamma$ the onto homomorphism given by $\pi(\bar{s})=s$, for every $s\in S$. 
\begin{lemma}\label{equivalence}
The group $\Gamma$ is {\it\emph{P}-stable} if and only if the following condition is satisfied:\\ $(\star)$
for every $T\subset S$ finite and $\varepsilon>0$, there are $E\subset\ker{\pi}$ finite and $\delta>0$ such that for any finite set $X$ and homomorphism $\rho:\mathbb F_S\rightarrow\emph{Sym}(X)$ satisfying $\emph{d}_{\emph{H}}(\rho(g),\emph{Id}_{X})\leq\delta$, for all $g\in E$, there is a homomorphism $\tau:\Gamma\rightarrow\emph{Sym}(X)$  satisfying $\emph{d}_{\emph{H}}(\rho(\bar{s}),\tau(s))\leq\varepsilon$, for all $s\in T$. 

Moreover, if $S$ is finite, then $\Gamma$ is \emph{P}-stable if and only if $(\star)$ is satisfied for $T=S$.
\end{lemma}

{\it Proof.} In the above notation, let $E_n\subset \ker{\pi}$ be an increasing sequence of sets with $\cup_nE_n=\ker(\pi)$. Let $p:\Gamma\rightarrow\mathbb F_S$ be a map such that $p(s)=\bar{s}$, for any $s\in S$, and $\pi(p(g))=g$, for any $g\in\Gamma$.

If $(\star)$ fails, then there exist $T\subset S$ finite, $\varepsilon>0$ and homomorphisms $\rho_n:\mathbb F_S\rightarrow\text{Sym}(X_n)$, with $X_n$ finite, such that $\max\{\text{d}_{\text{H}}(\rho_n(g),\text{Id}_{X_n})| g\in E_n\}\leq\frac{1}{n}$ and $\max\{\text{d}_{\text{H}}(\rho_n(\bar{s}),\tau_n(s))| s\in T\}>\varepsilon$, for any $n\in\mathbb N$ and homomorphism $\tau_n:\Gamma\rightarrow\text{Sym}(X_n)$. Define $\sigma_n:\Gamma\rightarrow \text{Sym}(X_n)$ by  $\sigma_n(g)=\rho_n(p(g))$. If $g,h\in\Gamma$, then $p(gh)^{-1}p(g)p(h)\in\ker{\pi}$, hence $p(gh)^{-1}p(g)p(h)\in E_{n_0}$, for some $n_0\in\mathbb N$. Therefore,  $$\text{$\text{d}_{\text{H}}(\sigma_n(gh),\sigma_n(g)\sigma_n(h))=\text{d}_{\text{H}}(\rho_n(p(gh)^{-1}p(g)p(h)),\text{Id}_{X_n})\leq\frac{1}{n}$, for every $n\geq n_0$.}$$ 
Then $(\sigma_n)_{n\in\mathbb N}$ is an asymptotic homomorphism of $\Gamma$.
On the other hand, as $\sigma_n(s)=\rho_n(\bar{s})$, for every $s\in S$, we get that $\max\{\text{d}_{\text{H}}(\sigma_n(s),\tau_n(s))|s\in T\}>\varepsilon$, for any homomorphism $\tau_n:\Gamma\rightarrow\text{Sym}(X_n)$ and $n\in\mathbb N$. This implies that $\Gamma$ is not P-stable.

Conversely, if $\Gamma$ is not P-stable, then there exist an asymptotic homomorphism
 $\sigma_n:\Gamma\rightarrow\text{Sym}(X_n)$, a finite set $T\subset S$ and $\varepsilon>0$ such that $\max_{s\in T}\text{d}_{\text{H}}(\sigma_n(s),\tau_n(s))>\varepsilon$, for any $n\in\mathbb N$ and  homomorphism $\tau_n:\Gamma\rightarrow\text{Sym}(X_n)$.
 Let $\rho_n:\mathbb F_S\rightarrow\text{Sym}(X_n)$ be the homomorphism given by $\rho_n(\bar{s})=\sigma_n(s)$, for all $s\in S$. Let $g\in\ker{\pi}$ and write $g=\bar{s}_1^{\varepsilon_1}...\bar{s}_k^{\varepsilon_k}$, for $s_1,...,s_k\in S$ and $\varepsilon_1,...,\varepsilon_k\in\{\pm 1\}$. Then $\rho_n(g)=\sigma_n(s_1)^{\varepsilon_1}...\sigma_n(s_k)^{\varepsilon_k}$. Since $s_1^{\varepsilon_1}...s_k^{\varepsilon_k}=e$ and $(\sigma_n)_{n\in\mathbb N}$ is an asymptotic homomorphism, we get that $\text{d}_{\text{H}}(\rho_n(g),\text{Id}_{X_n})\rightarrow 0$. Since $\max_{s\in T}\text{d}_{\text{H}}(\rho_n(\bar{s}),\tau_n(s))>\varepsilon$, for any $n\in\mathbb N$ and homomorphism $\tau_n:\Gamma\rightarrow\text{Sym}(X_n)$, we get that $(\star)$ is not satisfied. This finishes the proof of the lemma.
 \hfill$\blacksquare$

\subsection{P-stability vs. (very) flexible P-stability}
\begin{lemma}\label{vs}
Let $\Gamma$ be a countable group.
\begin{enumerate}
\item If $\Gamma$ is amenable, then it is \emph{P}-stable if and only if it is flexibly \emph{P}-stable.
\item If $\Gamma$ has property $(\tau)$, then it is flexibly \emph{P}-stable if and only if it is very flexibly \emph{P}-stable.
\end{enumerate}
\end{lemma}

{\it Proof.}
(1) Assume that $\Gamma$ is a flexibly P-stable amenable group.  In order to conclude that $\Gamma$ is P-stable, it is sufficient to prove the following claim:

{\bf Claim.} Let  $\sigma_n:\Gamma\rightarrow\text{Sym}(X_n)$ be an asymptotic homomorphism and $0<\varepsilon<1$. Then we can find a subsequence $(\sigma_{n_k})$ of $(\sigma_n)$ and homomorphisms $\tau_k:\Gamma\rightarrow\text{Sym}(X_{n_k})$, for any $k\in\mathbb N$, such that $\limsup_{k\rightarrow\infty}\text{d}_{\text{H}}(\sigma_{n_k}(g),\tau_k(g))\leq\varepsilon$, for every  $g\in\Gamma$.

To prove this claim we treat separately two cases.
Firstly, assume that $N:=\sup_n|X_n|<+\infty$. Since $\Gamma$ is flexibly P-stable, there are homomorphisms $\tau_n:\Gamma\rightarrow\text{Sym}(Y_n)$, with $Y_n\supset X_n$ finite, such that $|Y_n|/|X_n|\rightarrow 1$ and $\text{d}_{\text{H}}(\sigma_n(g),\tau_n(g)_{|X_n})\rightarrow 0$, for every $g\in\Gamma$. Thus, $|Y_n|/|X_n|<1+\frac{1}{N}$ and therefore $Y_n=X_n$, for $n$ large. This clearly implies the claim.

Secondly, assume that  $\sup_n|X_n|=+\infty$. After replacing $(\sigma_n)$ with a subsequence, we may suppose that $|X_n|\rightarrow +\infty$. 
Since $\Gamma$ is amenable, by using Ornstein and Weiss' theorem \cite{OW80} (similarly to the proof of \cite[Proposition 6.5]{BLT18}),  we can find a subsequence $(\sigma_{n_k})$ of $(\sigma_n)$ and $A_k\subset X_{n_k}$, for any $k\in\mathbb N$, such that $|\sigma_{n_k}(g)A_k\triangle A_k|/|X_{n_k}|\rightarrow 0$, for every $g\in\Gamma$, and $|A_{k}|/|X_{n_k}|\rightarrow \lambda:=1-\varepsilon$. 

For $k\in\mathbb N$, let $\rho_k:\Gamma\rightarrow\text{Sym}(A_k)$ be a map such that $\rho_k(g)$ agrees with $\sigma_{n_k}(g)$ on $A_k\cap\sigma_{n_k}(g)^{-1}A_k$, for every $g\in\Gamma$. Then $(\rho_k)$ is an asymptotic homomorphism. Since $\Gamma$ is flexibly P-stable, there are  $Y_k\supset A_k$ finite and homomorphisms $\zeta_k:\Gamma\rightarrow \text{Sym}(Y_k)$ such that $\text{d}_{\text{H}}(\rho_k(g),\zeta_k(g)_{|A_k})\rightarrow 0$, for every $g\in\Gamma$, and $|Y_k|/|A_k|\rightarrow 1$. Since $|A_k|/|X_{n_k}|\rightarrow \lambda<1$, we have $|Y_k|<|X_{n_k}|$ and so we may assume that $Y_k\subset X_{n_k}$, for $k$ large. If $\tau_k:\Gamma\rightarrow \text{Sym}(X_{n_k})$ is the homomorphism given by ${\tau_k(g)}_{|Y_k}=\zeta_k(g)$ and ${\tau_k(g)}_{|X_{n_k}\setminus Y_k}=\text{Id}_{X_{n_k}\setminus Y_k}$, then  $\limsup_{k\rightarrow\infty}\text{d}_{\text{H}}(\sigma_{n_k}(g),\tau_k(g))\leq\lim_{k\rightarrow\infty}|X_{n_k}\setminus A_k|/|X_{n_k}|=\varepsilon$, for every $g\in\Gamma$. This finishes the proof of the claim and of part (1).

(2) Assume that $\Gamma$ is a very flexibly P-stable group with property $(\tau)$. Let $\sigma_n:\Gamma\rightarrow\text{Sym}(X_n)$ be an asymptotic homomorphism. Then we can find homomorphisms $\tau_n:\Gamma\rightarrow\text{Sym}(Y_n)$, with $Y_n\supset X_n$ finite, such that $\text{d}_{\text{H}}(\sigma_n(g),\tau_n(g)_{|X_n})\rightarrow 0$, for any $g\in\Gamma$. 
Since $$\{x\in X_n\mid\tau_n(g)x\notin X_n\}\subset\{x\in X_n\mid\sigma_n(g)x\not=\tau_n(g)x\},$$ we get that $|\tau_n(g)X_n\triangle X_n|/|X_n|\rightarrow 0$, for any $g\in\Gamma$. Since $\Gamma$ has property $(\tau)$, Lemma \ref{almostinv} implies that $\sup_{g\in\Gamma}|\tau_n(g)X_n\triangle X_n|/|X_n|\rightarrow 0$. By Lemma \ref{component}, we can find a $\tau_n(\Gamma)$-invariant set $Z_n\subset Y_n$ such that $|Z_n\triangle X_n|/|X_n|\rightarrow 0$. Let $T_n=X_n\cup Z_n$ and $\rho_n:\Gamma\rightarrow \text{Sym}(T_n)$ be the homomorphism given by $\rho_n(g)_{|Z_n}=\tau_n(g)_{|Z_n}$ and $\rho_n(g)_{|X_n\setminus Z_n}=\text{Id}_{X_n\setminus Z_n}$. Then  we have $X_n\subset T_n$, $|T_n|/|X_n|\rightarrow 1$, and $\text{d}_{\text{H}}(\sigma_n(g),\rho_n(g)_{|X_n})\rightarrow 0$, for every $g\in\Gamma$. This shows that $\Gamma$ is flexibly P-stable.
\hfill$\blacksquare$

\subsection{Subgroups of finite index and very flexible P-stability} We end this section by proving that very flexible P-stability passes to subgroups  of finite index:

\begin{lemma}\label{finindex}
Let $\Gamma_0<\Gamma$ be a finite index inclusion of countable groups. 
If $\Gamma$ is very flexibly \emph{P}-stable, then so is $\Gamma_0$.
Moreover, if $\Gamma$ is weakly very flexibly \emph{P}-stable, then  so is $\Gamma_0$.
\end{lemma} 
{\it Proof.}  The proof is based on a simple induction argument (compare with \cite[Proposition 4.12]{ESS18}).
Assume that $\Gamma_0$ is not very flexibly \text{P}-stable. Then there exists an asymptotic homomorphism $\sigma_n:\Gamma_0\rightarrow\text{Sym}(X_n)$ a finite set $F\subset\Gamma_0$ and $\delta>0$ such that for any sequence of sets $Y_n\supset X_n$ and homomorphisms $\tau_n:\Gamma_0\rightarrow\text{Sym}(Y_n)$ we have that $\max\{\text{d}_{\text{H}}(\sigma_n(g),\tau_n(g)_{|X_n})\mid g\in F\}\geq\delta$, for all $n$.

Let $s:\Gamma/\Gamma_0\rightarrow\Gamma$ be a map such that $s(e\Gamma_0)=e$ and $s(g\Gamma_0)\in g\Gamma_0$, for all $g\in\Gamma$. Then $c:\Gamma\times\Gamma/\Gamma_0\rightarrow\Gamma_0$ given by $c(g,h\Gamma_0)=s(gh\Gamma_0)^{-1}g\;s(h\Gamma_0)$ is a cocycle for the left multiplication action $\Gamma\curvearrowright \Gamma/\Gamma_0$. For every $n$, we put $\widetilde X_n=\Gamma/\Gamma_0\times X_n$ and define $\widetilde\sigma_n:\Gamma\rightarrow\text{Sym}(\widetilde X_n)$ by letting $$\widetilde\sigma_n(g)(h\Gamma_0,x)=(gh\Gamma_0,\sigma_n(c(g,h\Gamma_0))x).$$
Then a direct computations shows that for every $g,h\in\Gamma$ we have $$\text{d}_{\text{H}}(\widetilde\sigma_n(gh),\widetilde\sigma_n(g)\widetilde\sigma_n(h))=\frac{1}{|\Gamma/\Gamma_0|}\;
\sum_{k\Gamma_0\in\Gamma/\Gamma_0}\text{d}_{\text{H}}(\sigma_n(c(gh,k\Gamma_0)),\sigma_n(c(g,hk\Gamma_0))\sigma_n(c(h,k\Gamma_0)).$$
Since $c(gh,k\Gamma_0)=c(g,hk\Gamma_0)c(h,k\Gamma_0)$, for all $k\Gamma_0\in\Gamma/\Gamma_0$, it follows that $\widetilde\sigma_n:\Gamma\rightarrow\text{Sym}(\widetilde X_n)$ is an asymptotic homomorphism.

Finally, consider a sequence of sets $Y_n\supset\widetilde X_n$ and homomorphisms $\tau_n:\Gamma\rightarrow\text{Sym}(Y_n)$. If  $g\in\Gamma_0$, then $\widetilde\sigma_n(g)$ leaves $e\Gamma_0\times X_n$ invariant and $\widetilde\sigma_n(g)(e\Gamma_0,x)=(e\Gamma_0,\sigma_n(g)x)$, for every $x\in X_n$. Thus, the restriction of ${\widetilde{\sigma_n}}_{|\Gamma_0}$ to $e\Gamma_0\times X_n$ can be identified to $\sigma_n$. Since ${\tau_n}_{|\Gamma_0}$ is a homomorphism, it follows that $\max\{\text{d}_{\text{H}}(\widetilde\sigma_n(g)_{|e\Gamma_0\times X_n},\tau_n(g)_{|e\Gamma_0\times X_n})\mid g\in F\}\geq\delta$. Thus,  $$\text{$\max\{\text{d}_{\text{H}}(\widetilde\sigma_n(g),\tau_n(g)_{|\widetilde X_n})\mid g\in F\}\geq \frac{\delta}{[\Gamma:\Gamma_0]}>0$, for all $n$},$$ which implies that $\Gamma$ is not very flexibly \text{P}-stable. This proves the main assertion. 

For the moreover assertion, assume the setting above and let $g\in\Gamma\setminus\{e\}$. Then we have that \begin{equation}\label{sigma_n}|\{\tilde x\in\widetilde X_n\mid\widetilde\sigma_n(g)\tilde x=\tilde x\}|=\sum_{h\Gamma_0\in\Gamma/\Gamma_0,gh\Gamma_0=h\Gamma_0}|\{x\in X_n\mid\sigma_n(c(g,h\Gamma_0))x=x\}|.\end{equation} 
If $h\Gamma_0\in\Gamma/\Gamma_0$ is such that $gh\Gamma_0=h\Gamma_0$, then we have $c(g,h\Gamma_0)=s(h\Gamma_0)^{-1}gs(h\Gamma_0)\not=e.$ 
Thus, if $\sigma_n:\Gamma_0\rightarrow\text{Sym}(X_n)$ is a sofic approximation of $\Gamma_0$, then using \eqref{sigma_n} it follows that $\widetilde\sigma_n:\Gamma\rightarrow\text{Sym}(\widetilde X_n)$ is a sofic approximation of $\Gamma$, and repeating the above argument implies the moreover assertion.
\hfill$\blacksquare$

\section{Permutation groups almost commuting with the regular representation}

The main goal of this section is to prove the following result. This implies that any group of permutations of a finite group $G$ that ``almost commutes" with the left regular representation of $G$ must arise from the right regular representation of $G$. More generally, we get precise structural information about any permutation group of a set containing $G$ whose restriction to $G$ almost commutes with the left regular representation of $G$. This generalization will be crucial later on in allowing us to prove that certain product groups are not very flexibly P-stable. 

\begin{theorem}\label{almost} Let $G$ be a finite group, $S$ be a set of generators and put $\kappa:=\kappa(G,S)$. Denote by $\alpha,\beta:G\rightarrow\emph{Sym}(X)$ the left and right multiplication actions of $G$ on $X:=G$. Let $Y$ be a finite set containing $X$ and $K<\emph{Sym}(Y)$ be a subgroup.
Let  $\varepsilon\in (0,\frac{\kappa^4}{200})$ and assume that
$$\text{$|\{x\in X\cap k^{-1}X\mid \alpha(g) kx\not=k\alpha(g)x\}|\leq\varepsilon\; |X|$, for all $g\in S$ and $k\in K$.}$$

Then $K_0=\{k\in K\mid |X\cap kX|\geq \frac{|X|}{2}\}$ is a subgroup of $K$. 

Moreover, we can find a homomorphism $\delta:K_0\rightarrow G$, a $K_0$-invariant set $X_1\subset Y$, a $\beta(\delta(K_0))$-invariant set $X_2\subset X$, and a bijection $\varphi:X_1\rightarrow X_2$ such that 
\begin{enumerate}
\item $|X\setminus X_1|=|X\setminus X_2|<\frac{4162}{\kappa^4}\;\varepsilon\;|X|$,
\item $|\{x\in X_1\mid \varphi(x)\not=x\}|\leq \frac{2048}{\kappa^4} \;\varepsilon \; |X|$, and
\item $\varphi\circ k_{|X_1}=\beta(\delta(k)))\circ\varphi$, for all $k\in K_0$. 
\end{enumerate}

\end{theorem}

The proof of Theorem \ref{almost} relies on the following two lemmas.

\begin{lemma}\label{commutant}\emph{\cite{Th10}} Let $G$ be a finite group, $S$ be a set of generators and put $\kappa:=\kappa(G,S)$. Denote by  $\alpha,\beta:G\rightarrow\emph{Sym}(G)$  the left and right multiplication actions of $G$ on itself. Then for every  $\varphi\in \emph{Sym}(G)$, there exists $h\in G$ such that $\kappa^2\cdot \emph{d}_{\emph{H}}(\varphi, \beta(h))\leq 4\;\max_{g\in S}\emph{d}_{\emph{H}}(\alpha(g)\circ\varphi,\varphi\circ\alpha(g))$. \end{lemma} 

After proving Lemma \ref{commutant}, we realized that it also follows from the proof of \cite[Theorem 2.2]{Th10}. Nevertheless, we include a self-contained proof for completeness. 

{\it Proof.}  Let  $\varphi\in\text{Sym}(G)$ and put $\varepsilon=\max_{g\in S}{\text d}_{\text H}(\alpha(g)\circ\varphi,\varphi\circ\alpha(g))$. Consider the unitary representation of $G$ on B$(\ell^2(G))$ given by $g\cdot T=\alpha(g)T\alpha(g)^*$, where we view $\text{Sym}(G)$ as a subgroup of $\text{U}(\ell^2(G))$.  Lemma \ref{almostinv} implies that $$\kappa\;\max_{g\in G}\|\alpha(g)\circ\varphi-\varphi\circ\alpha(g)\|_{\text{HS}}\leq 2\;\max_{g\in S}\|\alpha(g)\circ\varphi-\varphi\circ\alpha(g)\|_{\text{HS}}.$$ Recalling that $\|\sigma-\tau\|_{{\text HS}}=\sqrt{2\;{\text d}_{\text H}(\sigma,\tau)}$, for all $\sigma,\tau\in \text{Sym}(G)$, the last inequality rewrites as ${\text d}_{\text H}(\alpha(g)\circ\varphi,\varphi\circ\alpha(g))\leq \frac{4\varepsilon}{\kappa^2}$, for every $g\in G$. Equivalently, we have $|\{x\in G\mid\varphi(gx)\not= g\varphi(x)\}|\leq\frac{4\varepsilon}{\kappa^2}$, for every $g\in G$, and hence $$\sum_{x\in G}|\{g\in G\mid\varphi(gx)\not=g\varphi(x)\}|=\sum_{g\in G}|\{x\in G\mid\varphi(gx)\not=g\varphi(x)\}|\leq\frac{4\varepsilon}{\kappa^2}\;|G|.$$ Thus, there exists $x\in G$ such that $|\{g\in G\mid\varphi(gx)\not=g\varphi(x)\}|\leq\frac{4\varepsilon}{\kappa^2}$. Therefore, $h=\varphi(x)^{-1}x\in G$ satisfies ${\text d}_{\text H}(\varphi, \beta(h))\leq\frac{4\varepsilon}{\kappa^2}$. \hfill$\blacksquare$

\begin{lemma}\label{conjugacy}
Let $X$ be a finite set, $K$ a group and $\alpha_1,\alpha_2:K\rightarrow{\text Sym}(X)$ homomorphisms. Assume that $\emph{d}_{\emph{H}}(\alpha_1(k),\alpha_2(k))\leq\varepsilon$, for all $k\in K$, for some $\varepsilon>0$. 

Then there exist an $\alpha_1(K)$-invariant set $X_1\subset X$, an $\alpha_2(K)$-invariant set $X_2\subset X$, and a bijection $\varphi:X_1\rightarrow X_2$ such that
 $|X\setminus X_1|=|X\setminus X_2|\leq 16\varepsilon\; |X|$, 
$|\{x_1\in X_1\mid \varphi(x_1)\not=x_1\}|\leq 16\varepsilon\; |X|$, 
and
 $$\text{$\varphi\circ\alpha_1(k)_{|X_1}=\alpha_2(k)\circ\varphi$, for all $k\in K$}.$$
Moreover, if $\varepsilon<\frac{1}{16}$ and $\alpha_1$ is transitive, then $\alpha_1$ and $\alpha_2$ are conjugate.
\end{lemma}

{\it Proof.} We follow closely the proofs of \cite[Lemma 2.5]{Hj03} and \cite[Theorem 1.3]{Io06}.
We start by defining $V=\frac{1}{|K|}\sum_{k\in K}\alpha_2(k)^{-1}\circ\alpha_1(k)\in \text{B}(\ell^2(X))$. Then $\alpha_2(k)^{-1}V\alpha_1(k)=V$, for every $k\in K$. Thus, the matrix coefficients $V_{x_1,x_2}=\langle V\delta_{x_1},\delta_{x_2}\rangle$ satisfy \begin{equation}\label{equiv} \text{$V_{x_1,x_2}=V_{\alpha_1(k)x_1,\alpha_2(k)x_2}$, for all $x_1,x_2\in X$ and $k\in K$.}\end{equation}

Since $\|\alpha_1(k)^{-1}\circ\alpha_2(k)-\text{Id}\|_{\text HS}=\sqrt{2\;{\text d}_{\text H}(\alpha_1(k),\alpha_2(k))}\leq\sqrt{2\varepsilon}$, for every $k\in K$, we deduce that $\|V-\text{Id}\|_{\text HS}\leq\sqrt{2\varepsilon}$. Equivalently, we have \begin{equation}\label{V}\frac{1}{|X|}\Big(\sum_{x_1\in X}|V_{x_1,x_1}-1|^2+\sum_{x_1,x_2\in X,x_1\not=x_2}|V_{x_1,x_2}|^2\Big)\leq 2\varepsilon. \end{equation}

Let $A$ be the set of $x_1\in X$ for which there exists a unique $x_2=\varphi(x_1)\in X$ such that $|V_{x_1,x_2}|>\frac{1}{2}$. 
Then equation \eqref{equiv} implies that $A$ is $\alpha_1(K)$-invariant and \begin{equation}\label{equiv_2}\text{$\varphi(\alpha_1(k)x_1)=\alpha_2(k)\varphi(x_1)$, for all $x_1\in A$ and $k\in K$.} \end{equation}
Moreover, $A$  contains the set $X_0$ of $x_1\in X$ such that $|V_{x_1,x_1}-1|^2+\sum_{x_2\in X, x_2\not=x_1}|V_{x_1,x_2}|^2<\frac{1}{4}$. On the other hand, \eqref{V} implies that $\frac{|X\setminus X_0|}{4|X|}\leq 2\varepsilon$. Thus, $|X\setminus A|\leq |X\setminus X_0|\leq 8\varepsilon\; |X|$. Similarly, the set $B$ of $x_2\in X$ for which there is a unique $x_1\in X$ with $|V_{x_1,x_2}|>\frac{1}{2}$ satisfies $|X\setminus B|\leq 8\varepsilon\; |X|$.

Define $X_1=\{x_1\in A\mid \varphi(x_1)\in B\}$ and $X_2=\varphi(X_1)$. Then the restriction of $\varphi$ to $X_1$ is one-to-one.  Since $B$ is $\alpha_2(K)$-invariant, \eqref{equiv_2} gives that $X_1$ is $\alpha_1(K)$-invariant and $X_2$ is $\alpha_2(K)$-invariant. Since $\varphi(x_1)=x_1$, for all $x_1\in X_0$, we get that $X_0\cap B\subset X_1$. Thus, $|X\setminus X_1|\leq |X\setminus X_0|+|X\setminus B|\leq 16\varepsilon\; |X|$ and $|\{x_1\in X_1\mid\varphi(x_1)\not=x_1\}|\leq |X_1\setminus (X_0\cap B)|\leq |X\setminus (X_0\cap B)|\leq 16\varepsilon\; |X|$.

If $\varepsilon<\frac{1}{16}$, then $|X\setminus X_1|\leq 16\;\varepsilon\;|X|<|X|$, and thus $X_1$ is non-empty. Since $X_1$ is $\alpha_1(K)$-invariant, if $\alpha_1$ is transitive, we get that $X_1=X$ and the moreover assertion follows.
\hfill$\blacksquare$


\subsection*{Proof of Theorem \ref{almost}.} We will first show that $K_0$ is a subgroup of $K$. The proof of this assertion is inspired by the proof of \cite[Theorem 2.4]{GTD15}. If $g\in S$ and $k\in K$, then  \begin{align*}\alpha(g)(X\cap kX)\setminus (X\cap kX)&=\alpha(g)(\{x\in X\cap kX\mid\alpha(g)x\notin X\cap kX\}) \\&=\alpha(g)k(\{x\in X\cap k^{-1}X\mid \alpha(g)kx\notin X\cap kX\})\\&\subset \alpha(g)k(\{x\in X\cap k^{-1}X\mid \alpha(g)kx\not=k\alpha(g)x\}),\end{align*}
and thus $|\alpha(g)(X\cap kX)\setminus (X\cap kX)|\leq\varepsilon\;|X|$. 

Therefore, if $k\in K_0$, then for every $g\in S$ we have $|\alpha(g)(X\cap kX)\triangle (X\cap kX)|\leq 2\varepsilon\; |X|\leq 4\varepsilon \;|X\cap kX|$. By applying Lemma \ref{expansion} to $X\cap kX\subset X$ we deduce that 
$$ \kappa^2\;|X\cap kX|\; |X\setminus kX|\leq  \max_{g\in S}|\alpha(g)(X\cap kX)\triangle (X\cap kX)|\; |X|\leq 4\varepsilon \; |X\cap kX|\; |X|.$$

Hence, $|X\setminus  kX|\leq \frac{4\varepsilon}{\kappa^2}\; |X|$ and thus \begin{equation}\label{invariance}\text{$|X\triangle kX|\leq \frac{8\varepsilon}{\kappa^2}\;|X|$, for every $k\in K_0$.}\end{equation}If $k,k'\in K_0$, then $|X\triangle k'kX|\leq |X\triangle k'X|+|k'X\triangle k'kX|= |X\triangle k'X|+|X\triangle kX|\leq\frac{16\varepsilon}{\kappa^2}\;|X|$,  thus $|X\cap k'kX|\geq (1-\frac{8\varepsilon}{\kappa^2})|X|\geq |X|/2$ since  $\kappa\leq 2$ and hence $\varepsilon<\frac{\kappa^4}{200}<\frac{\kappa^2}{16}$. This shows that $kk'\in K_0$ and therefore $K_0$ is a subgroup of $K$.

Secondly, we will prove the existence of a map $\delta:K_0\rightarrow G$ such that \begin{equation}\label{unif}\text{$|\{x\in X\mid kx\not=\beta(\delta(k))x\}|\leq \frac{64\varepsilon}{\kappa^4}\; |X|$, for every $k\in K_0$}.\end{equation}
To see this, let $k\in K_0$. Let $\tilde k\in\text{Sym}(X)$ such that $\tilde kx=kx$, for every $x\in X\cap k^{-1}X$. If $g\in S$, then since $\tilde k\alpha(g)x=k\alpha(g)x$, for all $x\in X\cap\alpha(g)^{-1}k^{-1}X$, by using the hypothesis, we get that 
$$|\{x\in X\mid\alpha(g)\tilde kx\not=\tilde k\alpha(g)x\}|\leq \varepsilon\; |X|+|X\setminus (k^{-1}X\cap\alpha(g)^{-1}k^{-1}X)|\leq \varepsilon\; |X|+2\; |X\setminus k^{-1}X|.$$ In combination with \eqref{invariance}  this gives that
$$\text{$|\{x\in X\mid \alpha(g)\tilde kx\not=\tilde k\alpha(g)x\}|\leq (1+\frac{8}{\kappa^2})\varepsilon\; |X|$, for every $g\in S$.}$$

Now,
 Lemma \ref{commutant} gives $\delta(k)\in G$ such that $|\{x\in X\mid \tilde kx\not=\beta(\delta(k))x\}|\leq\frac{4}{\kappa^2}(1+\frac{8}{\kappa^2})\varepsilon\; |X|.$ Together with \eqref{invariance} we get that \begin{align*}|\{x\in X\mid kx\not=\beta(\delta(k))x\}|&\leq |\{x\in X\mid \tilde kx\not=\beta(\delta(k))x\}|+|X\setminus k^{-1}X|\\&\leq  \frac{4}{\kappa^2}(1+\frac{8}{\kappa^2})\varepsilon\; |X|+\frac{4}{\kappa^2}\varepsilon\; |X|.\end{align*} Since $\kappa\leq 2$, \eqref{unif} follows.
 
 Thirdly, we claim that $\delta:K_0\rightarrow G$ is a homomorphism. Denote $X_k=\{x\in X\mid kx=\beta(\delta(k))x\}$ for $k\in K$. Given $k',k\in K$, we have that $\beta(\delta(k'k))x=k'kx=k'\beta(\delta(k))x=\beta(\delta(k'))\beta(\delta(k))x$, for every $x\in X_{k'k}\cap X_k\cap \beta(\delta(k))^{-1}X_{k'}$. Thus, by using \eqref{unif} we get that $$|\{x\in X\mid\beta(\delta(k'k))x\not=\beta(\delta(k'))\beta(\delta(k))x\}|\leq\frac{192\varepsilon}{\kappa^4}\; |X|.$$
 Since $\varepsilon<\frac{\kappa^4}{200}$, we get that there exists $x\in X$ such that $\beta(\delta(k'k))x=\beta(\delta(k'))\beta(\delta(k))x$. Equivalently, $x\delta(k'k)^{-1}=x\delta(k)^{-1}\delta(k')^{-1}$, and thus $\delta(k'k)=\delta(k')\delta(k)$, which proves that $\delta$ is a homomorphism.

Finally, we will derive the rest of the conclusion by applying Lemma \ref{conjugacy}.
First, note that equation \eqref{invariance} together with Lemma \ref{component} provides a $K_0$-invariant set $X_0\subset Y$ such that $|X_0\triangle X|\leq\frac{16\varepsilon}{\kappa^2}\; |X|$. We put $Z=X_0\cup X$ 
and define homomorphisms $\alpha_1,\alpha_2:K_0\rightarrow\text{Sym}(Z)$ by letting for every $k\in K_0$ $$\text{$\alpha_1(k)_{|X_0}=k_{|X_0}, \alpha_1(k)_{|Z\setminus X_0}=\text{Id}_{Z\setminus X_0}$ and $\alpha_2(k)_{|X}=\beta(\delta(k))_{|X}$, $\alpha_2(k)_{|Z\setminus X}=\text{Id}_{Z\setminus X}.$}$$

Since $\{x\in Z\mid \alpha_1(k)x\not=\alpha_2(k)x\}\subset (X_0\triangle X)\cup\{x\in X_0\cap X\mid kx\not=\beta(\delta(k))x\}$, \eqref{unif} implies that
\begin{equation}
\text{$|\{x\in Z\mid\alpha_1(k)x\not=\alpha_2(k)x\}|\leq\frac{16\varepsilon}{\kappa^2} \;|X|+\frac{64\varepsilon}{\kappa^4}\; |X|\leq\frac{128\varepsilon}{\kappa^4}\; |X|$, for every $k\in K_0$.}
\end{equation}
By applying Lemma \ref{conjugacy}, we find an $\alpha_1(K_0)$-invariant set $Z_1\subset Z$, an $\alpha_2(K_0)$-invariant set $Z_2\subset Z$, and a bijection $\varphi:Z_1\rightarrow Z_2$ such that $|Z\setminus Z_1|\leq\frac{16\cdot 128\varepsilon}{\kappa^4}\; |X|, |\{z\in Z\mid\varphi(z)\not=z\}|\leq \frac{16\cdot 128\varepsilon}{\kappa^4}\;|X|$ and $\varphi\circ\alpha_1(k)_{|Z_1}=\alpha_2(k)\circ\varphi$, for all $k\in K_0$.
  Then $Z_1\cap X_0$ is  $\alpha_1(K_0)$-invariant and $Z_2\cap X$ is  $\alpha_2(K_0)$-invariant. Thus, $X_1=(Z_1\cap X_0)\cap\varphi^{-1}(Z_2\cap X)$  is $\alpha_1(K_0)$-invariant and $X_2=\varphi(X_1)$ is $\alpha_2(K_0)$-invariant. Since $X_1\subset X_0$ and $X_2\subset X$, we get that $X_1$ is $K_0$-invariant, $X_2$ is $\beta(\delta(K_0))$-invariant, and $\varphi\circ k_{|X_1}=\beta(\delta(k))\circ\varphi_{|X_1}$, for all $k\in K_0$. This proves condition {\it (3)} for $\varphi_{|X_1}$.
  
   In order to complete the proof, it remains to establish conditions {\it (1)} and {\it (2)}. First, we note that 
 $$|\{x\in X_1\mid\varphi(x)\not=x\}|\leq |\{z\in Z\mid\varphi(z)\not=z\}|\leq\frac{2048\varepsilon}{\kappa^4}\; |X|.$$
 Second, since $X_2=\varphi(Z_1\cap X_0)\cap (Z_2\cap X)$, we also get that \begin{align*}|X\setminus X_2|&\leq |Z\setminus Z_1|+|Z\setminus X_0|+|Z\setminus Z_2|+|Z\setminus X|
 \\&=2|Z\setminus Z_1|+|X_0\triangle X|
 \\ &\leq\frac{32\cdot 128\varepsilon}{\kappa^4}\; |X|+\frac{16\varepsilon}{\kappa^2}\; |X|\\&\leq\frac{4162\varepsilon}{\kappa^4}\; |X|,\end{align*}
 which finishes the proof.
\hfill$\blacksquare$

\section{A rigidity result for asymptotic homomorphisms}

In this section we prove the following consequence of Theorem \ref{almost}. 
\begin{theorem}\label{conjugation}
Let $\Gamma$ and $\Lambda$ be finitely generated groups. Assume that $\Gamma$ has property ($\tau$) with respect to a sequence of finite index normal subgroups $\{\Gamma_n\}_{n=1}^{\infty}$.
For every $n$, denote $X_n=\Gamma/\Gamma_n$, let $p_n:\Gamma\rightarrow X_n$ be the quotient homomorphism and $q_n:\Lambda\rightarrow X_n$ be a homomorphism.

Assume that $\sigma_n:(\Gamma*\mathbb Z)\times\Lambda\rightarrow\emph{Sym}(X_n)$, $n\in\mathbb N$, is an asymptotic homomorphism such that
\begin{enumerate}
\item For every $n\in\mathbb N$, we have $\sigma_n(g,h)x=p_n(g)xq_n(h)^{-1}$, for all $g\in\Gamma, h\in\Lambda, x\in X_n$.

\item For every $n\in\mathbb N$, there exist a finite set $Y_n$ which contains $X_n$ and a homomorphism $\tau_n:(\Gamma*\mathbb Z)\times\Lambda\rightarrow\emph{Sym}(Y_n)$ such that $\lim\limits_{n\rightarrow\infty}\emph{d}_{\emph{H}}(\sigma_n(g),\tau_n(g)_{|X_n})=0$, for all $g\in (\Gamma*\mathbb Z)\times\Lambda$.
\end{enumerate}

Then $\text{$\lim\limits_{n\rightarrow\infty}\Big(\max\{\emph{d}_{\emph{H}}(\sigma_n(t,e)\circ\sigma_n(e,h),\sigma_n(e,h)\circ\sigma_n(t,e))\mid h\in \Lambda\}\Big)=0$}$, for every $t\in\mathbb Z$.

Moreover, there exists a homomorphism $\sigma_n':(\Gamma*\mathbb Z)\times\Lambda\rightarrow\emph{Sym}(X_n)$ such that 
\begin{enumerate}[label=(\alph*)]
\item ${\sigma_n'}_{|\Gamma\times\Lambda}={\sigma_n}_{|\Gamma\times\Lambda}$, for every $n\in\mathbb N$, and
\item $\lim\limits_{n\rightarrow\infty}\emph{d}_{\emph{H}}(\sigma_n(g),\sigma_n'(g))=0$, for every $g\in (\Gamma*\mathbb Z)\times\Lambda$.
\end{enumerate}
\end{theorem}

{\it Proof.} Let $S$ and $T$ be finite sets of generators for $\Gamma$ and $\Lambda$, respectively. For $n\in\mathbb N$, we denote by $\beta_n:X_n\rightarrow\text{Sym}(X_n)$ the homomorphism given by $\beta_n(g)x=xp_n(g)^{-1}$. For ease of notation, we will write $g$ and $h$ instead of $(g,e)$ and $(e,h)$, for $g\in\Gamma*\mathbb Z$ and $h\in\Lambda$.

 In the first part of the proof we will use Theorem \ref{almost} to prove the following:

{\bf Claim.} For every $n$ large enough, there exist  a $\tau_n(\Lambda)$-invariant set $X_n'\subset Y_n$,  a subgroup $L_n<X_n$, a $\beta_n(L_n)$-invariant set $X_n''\subset X_n$ and 
 a bijection $\varphi_n:X_n'\rightarrow X_n''$ such that  
\begin{enumerate}
\item $\lim\limits_{n\rightarrow\infty}\frac{|X_n'|}{|X_n|}=\lim\limits_{n\rightarrow\infty}\frac{|X_n''|}{|X_n|}=1$,
\item $\lim\limits_{n\rightarrow\infty}\frac{1}{|X_n|}\; |\{x\in X_n'\mid\varphi_n(x)\not=x\}|=0$,
\item $\varphi_n\circ \tau_n(\Lambda)_{|X_n'}\circ\varphi_n^{-1}={\beta_n(L_n)}_{|X_n''}$.
\end{enumerate}

{\it Proof of the claim.}
For $n\geq 1$, we put $\varepsilon_n=2\;\max\{\text{d}_{\text{H}}(\sigma_n(g),\tau_n(g)_{|X_n})\mid g\in S\cup T\}$ and $$K_{n}=\{k\in \tau_n(\Lambda)\mid |X_n\cap kX_n|\geq\frac{|X_n|}{2}\}.$$

If $k\in\tau_n(\Lambda)$ and $g\in S$, then $k,\tau_n(g)\in\text{Sym}(Y_n)$ commute, thus $$\{x\in X_n\cap k^{-1}X_n\mid\sigma_n(g)kx\not=k\sigma_n(g)x\}\subset\{x\in X_n\cap k^{-1}X_n\mid \sigma_n(g)kx\not=\tau_n(g)kx\;\;\text{or}\;\;\sigma_n(g)x\not=\tau_n(g)x\}.$$
Therefore, for all $k\in\tau_n(\Lambda)$ and $g\in S$ we have \begin{equation}\label{S}\text{$|\{x\in X_n\cap k^{-1}X_n\mid\sigma_n(g)kx\not=k\sigma_n(g)x\}|\leq 2\;\text{d}_{\text{H}}(\sigma_n(g),\tau_n(g)_{|X_n})\;|X_n|\leq\varepsilon_n\;|X_n|$}.\end{equation}

Moreover, if $g\in T$, then $X_n\setminus\tau_n(g)^{-1}X_n=\{x\in X_n\mid\tau_n(g)x\not\in X_n\}\subset\{x\in X_n\mid\tau_n(g)x\not=\sigma_n(g)x\}$, and therefore
\begin{equation}\label{T} |X_n\cap\tau_n(g)X_n|=|X_n|-|X_n\setminus\tau_n(g)^{-1}X_n|\geq (1-\varepsilon_n)\; |X_n|.
\end{equation}

Since $\Gamma$ has property $(\tau)$ with respect to $\{\Gamma_n\}$ we have $\kappa:=\inf_n\kappa(X_n,p_n(S))>0.$ 
Since $\lim\limits_{n\rightarrow\infty}\varepsilon_n=0$, we have $\varepsilon_n<\min\{\frac{\kappa^4}{200},\frac{1}{2}\}$ for $n$ large enough.  By \eqref{S}, we can apply Theorem \ref{almost} to deduce that $K_{n}$ is a subgroup of $\tau_n(\Lambda)$ and there exist a $K_{n}$-invariant set $X_n'\subset Y_n$,  a subgroup $L_n<X_n$, a $\beta_n(L_n)$-invariant subset $X_n''\subset X_n$ and 
 a bijection $\varphi_n:X_n'\rightarrow X_n''$ such that  
\begin{itemize}
\item $|X_n'|=|X_n''|>(1-\frac{4162\;\varepsilon_n}{\kappa^4})|X_n|$,
\item $ |\{x\in X_n'\mid\varphi_n(x)\not=x\}|\leq \frac{2048\;\varepsilon_n}{\kappa^4}\;|X_n|$,
\item $\varphi_n\circ {K_{n}}_{|X_n'}\circ\varphi_n^{-1}={\beta_n(L_n)}_{|X_n''}$.
\end{itemize}
Since $\varepsilon_n<\frac{1}{2}$, \eqref{T} guarantees that $\tau_n(T)\subset K_n$. 
Since $K_n$ is a subgroup of $\tau_n(\Lambda)$ and $T$ generates $\Lambda$, we derive that $K_n=\tau_n(\Lambda)$. Since $\lim\limits_{n\rightarrow\infty}\varepsilon_n=0$, the claim follows. \hfill$\square$

Secondly, we claim that \begin{equation}\label{claimunu}\text{$\sigma_n(\Lambda)\subset \beta_n(L_n)$, for every $n$ large enough.}\end{equation}
To see this, let  $h\in T$. Then $\sigma_n(h)=\beta_n(q_n(h))$ and thus $\lim\limits_{n\rightarrow\infty}\text{d}_{\text{H}}(\beta_n(q_n(h)),\tau_n(h)_{|X_n})=0$. On the other hand, conditions (1)-(3) from above imply that we can find a sequence $h_n\in L_n$ such that $\lim\limits_{n\rightarrow\infty}\text{d}_{\text{H}}(\tau_n(h)_{|X_n},\beta_n(h_n))=0$.
Thus, we derive that $\lim\limits_{n\rightarrow\infty}\text{d}_{\text{H}}(\beta_n(q_n(h)),\beta_n(h_n))=0$. Since $\text{d}_{\text{H}}(\beta(k),\beta(k'))=\delta_{k,k'}$, for all $k,k'\in X_n$, we get that $q_n(h)=h_n\in L_n$, for large enough $n$. Since this holds for every $h\in T$, and $T$ is finite and generates $\Lambda$, the claim made in \eqref{claimunu} follows.

Thirdly, we claim that if $g\in\Gamma*\mathbb Z$, then $\sigma_n(g)$ asymptotically commutes with $\beta_n(L_n)$:
\begin{equation}\label{claimdoi} \lim\limits_{n\rightarrow\infty}\Big(\max\{\text{d}_{\text{H}}(\sigma_n(g)\circ\beta_n(h),\beta_n(h)\circ\sigma_n(g))\mid h\in L_n\}\Big)=0
\end{equation}
To see this, let $h_n\in L_n$, for every $n$. Condition (3) implies that $\beta_n(h_n)_{|X_n''}=\varphi_n\circ\tau_n(k_n)_{|X_n'}\circ\varphi_n^{-1}$, for some $k_n\in\Lambda$. By combining  (1) and (2) it follows that $\lim\limits_{n\rightarrow\infty}\text{d}_{\text{H}}(\beta_n(h_n),\tau_n(k_n)_{|X_n})=0$. On the other hand, we have $\lim\limits_{n\rightarrow\infty}\text{d}_{\text{H}}(\sigma_n(g),\tau_n(g)_{|X_n})=0$. Since $\tau_n(g)$ and $\tau_n(k_n)$ commute, we get that $\lim\limits_{n\rightarrow\infty}\text{d}_{\text{H}}(\sigma_n(g)\circ\beta_n(h_n),\beta_n(h_n)\circ\sigma_n(g))=0$. As this holds for any $h_n\in L_n$, claim \eqref{claimdoi} follows.

It is now clear that the combination of \eqref{claimunu} and \eqref{claimdoi} gives that \begin{equation}\label{claimtrei}\text{$ \lim\limits_{n\rightarrow\infty}\Big(\max\{\text{d}_{\text{H}}(\sigma_n(g)\circ\sigma_n(h),\sigma_n(h)\circ\sigma_n(g))\mid h\in\Lambda\}\Big)=0$, for every $g\in\Gamma*\mathbb Z$.}\end{equation}
Taking $g\in\mathbb Z$, this proves the main assertion. If $g\in\mathbb Z$ is a generator, then \eqref{claimtrei} together with Lemma \ref{commutant2} below implies the existence of $\sigma_n'(g)\in\text{Sym}(X_n)$ which commutes with $\sigma_n(\Lambda)$ such that $\lim\limits_{n\rightarrow\infty}\text{d}_{\text{H}}(\sigma_n'(g),\sigma_n(g))=0$. This implies the moreover assertion.
\hfill$\blacksquare$

In order to complete the proof of Theorem \ref{conjugation}, it remains to prove the following lemma.

\begin{lemma}\label{commutant2} Let $G$ be a finite group, $X$ a finite set, $\alpha:G\rightarrow\emph{Sym}(X)$ a homomorphism and $\varphi\in \emph{Sym}(X)$. Then there exists $\psi\in\emph{Sym}(X)$ which commutes with $\alpha(G)$ such that $$\emph{d}_{\emph{H}}(\varphi, \psi)\leq 32\;\max_{g\in G}\emph{d}_{\emph{H}}(\alpha(g)\circ\varphi,\varphi\circ\alpha(g)).$$

 \end{lemma}


{\it Proof of Lemma \ref{commutant2}.}  Put $\varepsilon=\max_{g\in G}{\text d}_{\text H}(\alpha(g)\circ\varphi,\varphi\circ\alpha(g))$. 
Then ${\text d}_{\text H}(\varphi^{-1}\circ\alpha(g)\circ\varphi,\alpha(g))\leq \varepsilon$, for any $g\in G$.  By applying Lemma \ref{conjugacy} to the homomorphisms $\varphi^{-1}\circ\alpha\circ\varphi,\alpha:G\rightarrow\text{Sym}(X)$ we find an $\alpha(G)$-invariant set $X_1\subset X$, an $\varphi^{-1}\alpha(G)\varphi$-invariant set $X_2\subset X$ and a bijection $\sigma:X_1\rightarrow X_2$ such that $|X\setminus X_1|\leq 16\varepsilon\; |X|$, $|\{x\in X_1\mid\sigma(x)\not=x\}|\leq 16\varepsilon\; |X|$ and $$\text{$\varphi^{-1}\circ\alpha(g)\circ\varphi\circ\sigma=\sigma\circ\alpha(g)_{|X_1}$, for all $g\in G$. }$$
Thus,  $X_3=\varphi(X_2)$ is $\alpha(G)$-invariant and the bijection $\tau=\varphi\circ\sigma:X_1\rightarrow X_3$ satisfies 
\begin{equation}\label{conjugat} \text{$\alpha(g)\circ\tau=\tau\circ\alpha(g)_{|X_1}$, for every $g\in G$.} \end{equation} 
Next, we say that two actions $\beta:G\rightarrow\text{Sym}(Y)$ and $\gamma:G\rightarrow\text{Sym}(Z)$ are conjugate if there exists a bijection $\rho:Y\rightarrow Z$ such that $\rho\circ\beta(g)=\gamma(g)\circ\rho$, for every $g\in G$. Let ${\text{Sub}_{\sim}(G)}$ be the set of equivalence classes of subgroups of $G$ modulo inner conjugacy. 
 For a subgroup $H<G$, denote by $\zeta({\beta})([H])$ the number of disjoint $\beta(G)$-orbits  $\beta(G)y$, with $y\in Y$, such that the restriction of $\beta$ to $\beta(G)y$ is conjugate to the action $G\curvearrowright G/H$. Then the conjugacy class of an action $\beta:G\rightarrow\text{Sym}(Y)$ is completely determined by the map $\zeta({\beta}):{\text{Sub}_{\sim}(G)}\rightarrow\mathbb N$. 

Finally, \eqref{conjugat} implies that  the restrictions of $\alpha$ to $X_1$ and $X_3$ are conjugate, hence $\zeta(\alpha_{|X_1})= \zeta(\alpha_{|X_3})$. This implies that  $\zeta(\alpha_{|X\setminus X_1})= \zeta(\alpha_{|X\setminus X_3})$ and so restrictions of $\alpha$ to $X\setminus X_1$ and $X\setminus X_3$ are conjugate.  In combination with \ref{conjugat}, we derive that there exists $\psi\in\text{Sym}(X)$ which commutes with $\alpha(G)$ (i.e., a self-conjugacy of $\alpha$) such that $\psi_{|X_1}=\tau$. Hence, $$|\{x\in X\mid\psi(x)\not=\varphi(x)\}|\leq |X\setminus X_1|+|\{x\in X_1\mid\sigma(x)\not=x\}|\leq 32\varepsilon\;|X|$$ and the conclusion follows. \hfill$\blacksquare$

\section{Construction of asymptotic homomorphisms}

This section is devoted to the construction of asymptotic homomorphisms. In the next section, we will combine this construction with Theorem \ref{conjugation} to deduce our main results.

\begin{lemma}
\label{tech2}
Let $\Gamma$ and $\Lambda$ be finitely generated groups. 
Let $\{\Gamma_n\}_{n=1}^{\infty}$ be a sequence of finite index normal subgroups of $\Gamma$,
put $X_n=\Gamma/\Gamma_n$ and denote by $p_n:\Gamma\rightarrow X_n$ the quotient homomorphism.
Assume that there exists a sequence of homomorphisms $q_n:\Lambda\rightarrow X_n$ such that $\Lambda$ does not have property $(\tau)$ with respect to the sequence $\{\ker(q_n)\}_{n=1}^{\infty}$. Let $t=\pm 1$ be a generator of $\mathbb Z$.

Then there exists an asymptotic homomorphism $\sigma_n:(\Gamma*\mathbb Z)\times\Lambda\rightarrow\emph{Sym}(X_n)$  such that
\begin{enumerate}
\item $\sigma_n(g,h)x=p_n(g)xq_n(h)^{-1}$, for all $g\in\Gamma,h\in\Lambda,x\in X_n$.
\item $\max\{\emph{d}_{\emph{H}}(\sigma_n(t,e)\circ\sigma_n(e,h),\sigma_n(e,h)\circ\sigma_n(t,e))\mid h\in\Lambda\}\geq\frac{1}{126}$, for infinitely many $n\geq 1$.
\end{enumerate}
\end{lemma}

{\it Proof.}
{The first part of the proof is devoted to the construction of  $\sigma_n$.} Let $\sigma_n:\Gamma\times\Lambda\rightarrow\text{Sym}(X_n)$ be  given by {\it (1)}. In order to extend $\sigma_n$ to an asymptotic homomorphism of $(\Gamma*\mathbb Z)\times\Lambda$ we will define $\sigma_n(t,e)\in\text{Sym}(X_n)$ such that $\lim\limits_{n\rightarrow\infty}\text{d}_{\text{H}}(\sigma_n(t,e)\circ\sigma_n(e,h),\sigma_n(e,h)\circ\sigma_n(t,e))=0$, for any $h\in\Lambda$.

To this end,  let $T\subset \Lambda$ be a finite generating set. Since $\Lambda$ does not have property $(\tau)$ with respect to $\{\ker(q_n)\}_{n=1}^{\infty}$ we have that $\inf_n\kappa(q_n(\Lambda),q_n(T))=0$. Thus, after passing to a subsequence, we may assume that $\lim\limits_{n\rightarrow\infty}\kappa(q_n(\Lambda),q_n(T))=0$. 
 Lemma \ref{AE} then implies that for every $n$ large enough there exists a set $C_n\subset q_n(\Lambda)$ such that \begin{equation}\label{C_n}\text{$\frac{1}{7}\leq\frac{|C_n|}{|q_n(\Lambda)|}\leq\frac{1}{6}$ and $\lim\limits_{n\rightarrow\infty}\frac{|C_nq_n(h)\triangle C_n|}{|q_n(\Lambda)|}=0$, for every $h\in\Lambda$.}\end{equation}

Let $Z_n\subset X_n$ be a set of representatives for the left cosets of $q_n(\Lambda)$. We define $B_n=Z_n\cdot C_n\subset X_n$, and claim that $B_n$ satisfies the following:
\begin{enumerate}[label=(\alph*)]
\item\label{a} $\frac{1}{7}\leq\frac{|B_n|}{|X_n|}\leq\frac{1}{6}$,  
\item $\lim\limits_{n\rightarrow\infty}\frac{|B_nq_n(h)\triangle B_n|}{|X_n|}=0$, for every $h\in\Lambda$, and
\item\label{c} $\frac{1}{|q_n(\Lambda)|}\;\sum_{h\in q_n(\Lambda)}|B_nh\cap Y|=\frac{|B_n|\; |Y|}{|X_n|}$, for every $Y\subset X_n$.\end{enumerate}

Indeed, (a) and (b) follow from \eqref{C_n}. To verify (c), note that if $x\in X_n$, then there is a unique $z\in Z_n$ such that $x^{-1}z\in q_n(\Lambda)$ and thus we have that $$|x^{-1}B_n\cap q_n(\Lambda)|=|\{h\in q_n(\Lambda)\mid h\in x^{-1}Z_n\cdot C_n\}|=|\{h\in q_n(\Lambda)\mid h\in (x^{-1}z)C_n\}|=|C_n|.$$
Therefore, we deduce that $$\sum_{h\in q_n(\Lambda)}|B_nh\cap Y|=\sum_{h\in q_n(\Lambda), x\in X_n}{\bf 1}_{B_nh}(x){\bf 1}_Y(x)=\sum_{x\in Y}|x^{-1}B_n\cap q_n(\Lambda)|=|C_n|\; |Y|$$
Since $|B_n|=\frac{|X_n|}{|q_n(\Lambda)|}\;|C_n|$, condition (c) is also satisfied.

Let $n$ large enough. Since $\sum_{g\in X_n}|B_n\setminus g^{-1}B_n|=|B_n|\cdot (|X_n|-|B_n|)$,  we can find $g_n\in X_n$ such that $A_n=B_n\setminus g_n^{-1}B_n$ satisfies $\frac{|A_n|}{|X_n|}\geq \frac{|X_n|-|B_n|}{|X_n|}\; \frac{|B_n|}{|X_n|}\geq\frac{5}{6}\cdot \frac{1}{7}=\frac{5}{42}.$ Moreover, $A_n\cap g_nA_n=\emptyset$ and
\begin{equation}\label{A_n}
\text{$\lim\limits_{n\rightarrow\infty}\frac{|A_nq_n(h)\triangle A_n|}{|X_n|}=0$, for every $h\in\Lambda$.}
\end{equation}

We are now ready to define $\sigma_n(t,e)\in\text{Sym}(X_n)$ by letting $$\sigma_n(t,e)x=\begin{cases}\text{$g_nx$, if $x\in A_n$,}\\ \text{$g_n^{-1}x$, if $x\in g_nA_n$,} \\ \text{$x$, otherwise.}  \end{cases}$$
Then for every $h\in\Lambda$ we have that $$(\sigma_n(t,e)\circ\sigma_n(e,h))x=\begin{cases}g_nxq_n(h)^{-1}, \text{if $x\in A_nq_n(h)$}, \\ g_n^{-1}xq_n(h)^{-1}, \text{if $x\in g_nA_nq_n(h)$}, \\ \text{$xq_n(h)^{-1}$, otherwise}\end{cases}$$
and $$(\sigma_n(e,h)\circ\sigma_n(t,e))x=\begin{cases}
\text{$g_nxq_n(h)^{-1}$, if $x\in A_n$,}\\ 
\text{$g_n^{-1}xq_n(h)^{-1}$, if $x\in g_nA_n$,}\\
\text{$xq_n(h)^{-1}$, otherwise.}
\end{cases} $$
These formulae easily imply that $\text{d}_{\text{H}}(\sigma_n(t,e)\circ\sigma_n(e,h),\sigma_n(e,h)\circ\sigma_n(t,e))$ is equal to \begin{equation}\label{commute}\begin{cases}\text{$\frac{2|A_n\setminus A_nq_n(h)|+|(A_n\cup g_nA_n)\setminus (A_n\cup g_nA_n)q_n(h)|}{|X_n|}$, if $g_n^2\not=e$}\\ \text{$\frac{2|(A_n\cup g_nA_n)\setminus (A_n\cup g_nA_n)q_n(h)|}{|X_n|}$, if $g_n^2=e$.}\end{cases}\end{equation}
Since $(A_n\cup g_nA_n)\setminus (A_n\cup g_nA_n)q_n(h)\subset (A_n\setminus A_nq_n(h))\cup g_n(A_n\setminus A_nq_n(h))$,  \eqref{A_n} implies that for all $h\in\Lambda$, $\lim\limits_{n\rightarrow\infty}\text{d}_{\text{H}}(\sigma_n(t,e)\circ\sigma_n(e,h),\sigma_n(e,h)\circ\sigma_n(t,e))=0$. This ends the first part of the proof.

 In the second part of the proof we will prove condition {\it (2)} from the conclusion. Let $n$ large enough. By using  \eqref{A_n},\eqref{commute}, that $\frac{|A_n|}{|X_n|}\geq\frac{5}{42}$ and that $A_n\subset B_n$, for all $h\in\Lambda$, we get that 
\begin{align}\begin{split}\label{d_H}&\text{d}_{\text{H}}(\sigma_n(t,e)\circ\sigma_n(e,h),\sigma_n(e,h)\circ\sigma_n(t,e))\\
&\geq \frac{|(A_n\cup g_nA_n)\setminus (A_n\cup g_nA_n)q_n(h)|}{|X_n|}\\&\geq\frac{|A_n|}{|X_n|}-\frac{|(A_n\cup g_nA_n)\cap (A_nq_n(h)\cup g_nA_nq_n(h))|}{|X_n|}\\&\geq\frac{5}{42}-\frac{|(B_n\cup g_nB_n)\cap (B_nq_n(h)\cup g_nB_nq_n(h))|}{|X_n|}\end{split}
\end{align}
Since $|(B_n\cup g_nB_n)\cap (B_nh\cup g_nB_nh)|\leq 2|B_nh\cap B_n|+|B_nh\cap g_nB_n|+|B_nh\cap g_n^{-1}B_n|$, by using condition \ref{c} we derive that $$\frac{1}{|q_n(\Lambda)|}\;\sum_{h\in q_n(\Lambda)}|(B_n\cup g_nB_n)\cap (B_nh\cup g_nB_nh)|\leq 4\;\frac{|B_n|^2}{|X_n|}.$$
Thus, there exists $h_n\in\Lambda$ such that $|(B_n\cup g_nB_n)\cap (B_nq_n(h_n)\cup g_nB_nq_n(h_n))|\leq 4\;\frac{|B_n|^2}{|X_n|}$. By combining this with \eqref{d_H} and the inequality $\frac{|B_n|}{|X_n|}\leq\frac{1}{6}$ from \ref{a}, it follows  that $h_n$ satisfies \begin{equation}
\text{d}_{\text{H}}(\sigma_n(t,e)\circ\sigma_n(e,h_n),\sigma_n(e,h_n)\circ\sigma_n(t,e)) \geq \frac{5}{42}-4\;\frac{|B_n|^2}{|X_n|^2}\geq\frac{5}{42}-\frac{4}{36}=\frac{1}{126}
\end{equation}
This proves condition {\it (2)} and finishes the proof.
\hfill$\blacksquare$

\section{Proofs of main results}

The proof of Theorem \ref{A} relies on the following result which puts together Theorems \ref{conjugation} and \ref{tech2}.

\begin{theorem}\label{tech}
Let $\Gamma$ and $\Lambda$ be finitely generated groups. Assume that
$\Gamma$ has  property $(\tau)$ with respect to a sequence $\{\Gamma_n\}_{n=1}^{\infty}$ of finite index normal subgroups. Suppose that there exist homomorphisms $q_n:\Lambda\rightarrow X_n$ such that $\Lambda$ does not have property $(\tau)$ with respect to the sequence $\{\ker(q_n)\}_{n=1}^{\infty}$.

Then $\Sigma\times\Lambda$ is not very flexibly \emph{P}-stable, for any finitely generated group $\Sigma$  which factors onto $\Gamma*\mathbb Z$.
\end{theorem}

{\it Proof.} Assume by contradiction that $\Sigma\times\Lambda$ is very flexibly P-stable. Let $\pi:\Sigma\rightarrow\Gamma*\mathbb Z$ be an onto homomorphism, and denote still by $\pi$ the product homomorphism $\pi\times\text{Id}_{\Lambda}:\Sigma\times\Lambda\rightarrow (\Gamma*\mathbb Z)\times\Lambda$. Let $t=\pm 1$ be a generator of $\mathbb Z$. Denote $X_n=\Gamma/\Gamma_n$ and let $p_n:\Gamma\rightarrow X_n$ be the quotient homomorphism. By Lemma \ref{tech2}, there exists an asymptotic homomorphism $\sigma_n:(\Gamma*\mathbb Z)\times\Lambda\rightarrow\text{Sym}(X_n)$  such that
\begin{enumerate}
\item $\sigma_n(g,h)x=p_n(g)xq_n(h)^{-1}$, for all $g\in\Sigma,h\in\Lambda,x\in X_n$.
\item $\max\{\text{d}_{\text{H}}(\sigma_n(t,e)\circ\sigma_n(e,h),\sigma_n(e,h)\circ\sigma_n(t,e))\mid h\in\Lambda\}\geq\frac{1}{126}$, for infinitely many $n\geq 1$.
\end{enumerate}

Then $\sigma_n\circ\pi:\Sigma\times\Lambda\rightarrow\text{Sym}(X_n)$ is an asymptotic homomorphism. Thus, since $\Sigma\times\Lambda$ is assumed very flexibly P-stable, for every $n\in\mathbb N$, we can find a finite set $Y_n\supset X_n$ and a homomorphism $\tau_n:\Sigma\times\Lambda\rightarrow\text{Sym}(Y_n)$ such that $\lim\limits_{n\rightarrow\infty}\text{d}_{\text{H}}(\sigma_n(\pi(g)),\tau_n(g)_{|X_n})=0$, for every $g\in\Sigma\times\Lambda$.

Since $\Gamma$ is finitely generated and $\pi$ is onto, we can find a finitely generated subgroup $\Delta<\Sigma$ and  $\tilde t\in\Sigma$ such that $\pi(\Delta)=\Gamma$ and $\pi(\tilde t)=t$. Let $\rho:\Delta*\mathbb Z\rightarrow \Sigma$ the homomorphism given by $\rho_{|\Delta}=\text{Id}_{\Delta}$ and $\rho(t)=\tilde t$. Denote still by $\rho$ the product homomorphism $\rho\times\text{Id}_{\Lambda}:(\Delta*\mathbb Z)\times\Lambda\rightarrow\Sigma\times\Lambda$. Then $\alpha_n:=\sigma_n\circ\pi\circ\rho:(\Delta*\mathbb Z)\times\Lambda\rightarrow\text{Sym}(X_n)$ is an asymptotic homomorphism which satisfies that $\lim\limits_{n\rightarrow\infty}\text{d}_{\text{H}}(\alpha_n(g),\tau_n(\rho(g))_{|X_n})=0$, for every $g\in(\Delta*\mathbb Z)\times\Lambda$.

 Now, note that (1) gives that $\alpha_n(g,h)x=(p_n\circ\pi)(g)xq_n(h)^{-1}$, for all $g\in\Delta,h\in\Lambda$  and $x\in X_n$.
Since $\Gamma$ has property $(\tau)$ with respect to $\{\Gamma_n\}_{n=1}^{\infty}$, $\Delta$ has property $(\tau)$ with respect to $\{\ker(p_n\circ\pi)\}_{n=1}^{\infty}$. Since $p_n\circ\pi:\Delta\rightarrow X_n$ is an onto homomorphism and $\tau_n\circ\rho:(\Delta*\mathbb Z)\times\Lambda\rightarrow \text{Sym}(Y_n)$ is a homomorphism, for all $n\in\mathbb N$, applying Theorem \ref{conjugation} 
to $\alpha_n:(\Delta*\mathbb Z)\times\Lambda\rightarrow\text{Sym}(X_n)$ gives that 
$$\lim\limits_{n\rightarrow\infty}\big(\max\{\text{d}_{\text{H}}(\alpha_n(t,e)\circ\alpha_n(e,h),\alpha_n(e,h)\circ\alpha_n(t,e))\mid h\in\Lambda\}\big)=0.$$
However, since $\alpha_n(t,e)=\sigma_n(t,e)$ and $\alpha_n(e,h)=\sigma_n(e,h)$, for every $h\in\Lambda$, this contradicts (2).
\hfill$\blacksquare$

\subsection*{Proof of Theorem \ref{A}} Let $\Sigma$ and $\Lambda$ be finitely generated groups such that $\Sigma$ admits a non-abelian free quotient and $\Lambda$ does not have property $(\tau)$.   Our goal is to prove that $\Sigma\times\Lambda$ is not very flexibly \text{P}-stable. By Lemma \ref{finindex} it suffices to find a finite index subgroup $\Sigma_0<\Sigma$ such that $\Sigma_0\times\Lambda$ is not very flexibly P-stable.

Let us first prove the conclusion in the case when $\Lambda$ admits an infinite cyclic quotient, since this requires less technology than the general case. Let $\rho:\Lambda\rightarrow\mathbb Z$ be an onto homomorphism. Since $\Sigma$ factors onto $\mathbb F_2$, it has a finite index subgroup $\Sigma_0$ which factors onto $\mathbb F_3$. Towards showing that $\Sigma_0\times\Lambda$ is not very flexibly P-stable, recall that $\Gamma=\mathbb F_2$  can be realized as a finite index subgroup of SL$_2(\mathbb Z)$, by letting for instance $\Gamma=\langle\begin{pmatrix}1&2\\ 0&1 \end{pmatrix},\begin{pmatrix}1&0\\2&1\end{pmatrix}\rangle$. Since SL$_2(\mathbb Z)$ has the Selberg property \cite{LW93} (i.e., property $(\tau)$ with respect to its congruence subgroups),  $\Gamma$ has property $(\tau)$ with respect to $\{\Gamma_n\}_{n=1}^{\infty}$, where $\Gamma_n=\Gamma\cap \big(\ker(\text{SL}_2(\mathbb Z)\rightarrow\text{SL}_2(\frac{\mathbb Z}{n\mathbb Z})\big)$. Let $p_n:\Gamma\rightarrow\Gamma/\Gamma_n$ be the quotient homomorphism. Let $\eta:\mathbb Z\rightarrow\Gamma$ the homomorphism given by $\eta(1)=\begin{pmatrix}1&2\\0&1 \end{pmatrix}$ and denote $q_n=p_n\circ\eta\circ\rho:\Lambda\rightarrow\Gamma/\Gamma_n$. Since $q_n$ factors through $\rho:\Lambda\rightarrow\mathbb Z$, for every $n$, and $\lim\limits_{n\rightarrow\infty}|q_n(\Lambda)|=+\infty$, it follows that $\Lambda$ does not have property $(\tau)$ with respect to $\{\ker(q_n)\}_{n=1}^{\infty}$. Since $\Sigma_0$ factors onto $\mathbb F_3=\Gamma*\mathbb Z$, Theorem \ref{tech} implies that $\Sigma_0\times\Lambda$ is not very flexibly P-stable.

In order to establish the general case we will use a theorem of Kassabov \cite[Theorem 2]{Ka05} which provides an integer $L\geq 2$ and onto homomorphisms $\pi_n:\mathbb F_L\rightarrow\text{Sym}(n)$, for every $n\in\mathbb N$, such that $\inf_n\kappa(\text{Sym}(n),\pi_n(S))>0$, where $S\subset\mathbb F_L$ is a free generating set. In other words, $\Gamma=\mathbb F_L$ has property $(\tau)$ with respect to $\{\ker(\pi_n)\}_{n=1}^{\infty}$. Since $\Sigma$ factors onto $\mathbb F_2$, it has a finite index subgroup $\Sigma_0$ which factors onto $\mathbb F_{L+1}$. We will show that $\Sigma_0\times\Lambda$ is not very flexibly P-stable. 

To this end, note that since $\Lambda$ does not have property $(\tau)$, there exists a sequence $\{\Lambda_n\}_{n=1}^{\infty}$ of finite index normal subgroups such that  $\lim\limits_{n\rightarrow\infty}\kappa(\Lambda/\Lambda_n,\delta_n(T))=0$, where $T\subset\Lambda$ is a finite generating set and $\delta_n:\Lambda\rightarrow\Lambda/\Lambda_n$ denote the quotient homomorphisms. For every $n\in\mathbb N$, put $G_n=\text{Sym}(\Lambda/\Lambda_n)$ and let $i_n:\Lambda/\Lambda_n\rightarrow G_n$ be the embedding given by left multiplication action of $\Lambda/\Lambda_n$ on itself. We denote $q_n=i_n\circ\delta_n:\Lambda\rightarrow G_n$.
Finally, we put $k_n=|\Lambda/\Lambda_n|$ and let $p_n:\Gamma\rightarrow G_n$ be the onto homomorphism obtained by composing $\pi_{k_n}:\Gamma\rightarrow\text{Sym}(n_k)$ with an isomorphism $\text{Sym}(k_n)\cong G_n$. By construction,  $\Gamma$ has property $(\tau)$ with respect to $\{\ker(p_n)\}_{n=1}$, while $\Lambda$ does not have property $(\tau)$ with respect to $\{\ker(q_n)\}_{n=1}^{\infty}$ (as $\ker(p_n)=\ker(\pi_{k_n})$ and $\ker(q_n)=\Lambda_n$, for every $n\in\mathbb N$). Since $\Sigma_0$ factors onto $\mathbb F_{L+1}=\Gamma*\mathbb Z$, Theorem \ref{tech} implies that $\Sigma_0\times\Lambda$ is not very flexibly stable. \hfill$\blacksquare$

\subsection*{Proof of Corollary \ref{B}} Since $\mathbb Z^d$ and $\mathbb F_n$ do not have property ($\tau$) for any integers $d,n\geq 1$, parts (1) and (2) follow  from Theorem \ref{A}. Let $m,n$ be integers such that $|m|=|n|\geq 2$. Then the Baumslag-Solitar group BS$(m,n)=\langle a,t|ta^mt^{-1}=a^n\rangle$ has a finite index subgroup isomorphic to $\mathbb F_{k}\times\mathbb Z$ for some $k\geq 2$ (see., e.g., \cite[Proposition 2.6]{Le05}). Since $\mathbb F_{k}\times\mathbb Z$ is not very flexbily P-stable by part (1), the same is true for BS$(m,n)$ by Lemma \ref{finindex}. This proves part (3). To prove part (4), let $n\geq 3$ be an integer. Recall that the pure braid group PB$_n$ has infinite center, Z(PB$_n)\cong\mathbb Z$, and admits a non-trivial splitting $\text{PB}_n\cong \text{PB}_n/\text{Z(PB}_n)\times\text{Z(PB}_n)$ (see \cite[Chapter 9]{FM11}). Since $\text{PB}_m$ factors onto $\text{PB}_{m-1}$, for any $m\geq 3$, and $\text{PB}_3\cong\mathbb F_2\times\mathbb Z$, we get that $\text{PB}_n$  factors onto $\mathbb F_2$. Thus, $\text{PB}_n/\text{Z(PB}_n)$ factors onto $\mathbb F_2$. Applying Theorem \ref{A} implies that $\text{PB}_n$ is not very flexibly P-stable. Since $\text{PB}_n$ is a finite index subgroup of $\text{B}_n$, the same holds for $\text{B}_n$ by Lemma \ref{finindex}.
\hfill$\blacksquare$

\subsection*{Proof of Theorem \ref{C}}
By the moreover assertion of Lemma \ref{finindex}, it suffices to prove that $\Sigma\times\Lambda$ is not weakly very flexibly P-stable, where $\Sigma=\mathbb F_m$ and $\Lambda=\mathbb Z^d$ or $\Lambda=\mathbb F_k$, for $m,k\geq 2$ and $d\geq 1$. Since any subgroup of index $2$ of $\mathbb F_m$ is isomorphic to $\mathbb F_{2m-1}$, using Lemma \ref{finindex} again we may assume that $m\geq 3$.  Let $\Gamma=\mathbb F_{m-1}$, so that $\Sigma=\Gamma*\mathbb Z$.
We view $\Sigma$ as a  finite index subgroup of $\text{SL}_2(\mathbb Z)$, and denote by $\pi_r:\text{SL}_2(\mathbb Z)\rightarrow\text{SL}_2(\mathbb Z/r\mathbb Z)$ the quotient homomorphism, for any prime $r$.

 We continue by treating two cases:

{\bf Case 1.} $\Lambda=\mathbb Z^d$, for some $d\geq 1$.

Fix $n\in\mathbb N$ and let $r_{n,0}, r_{n,1},..,r_{n,d}$ be $d+1$ distinct primes greater than $n$. Define $X_n=\prod_{i=1}^d\text{SL}_2(\mathbb Z/r_{n,i}\mathbb Z)$ and homomorphisms $p_n:\Gamma\rightarrow X_n, q_n:\Lambda\rightarrow X_n$ by letting for $g\in\Gamma$ and $(h_1,...,h_d)\in\Lambda$  $$\text{$p_n(g)=(\pi_{r_{n,1}}(g),...,\pi_{r_{n,d}}(g))$ and $q_n(h_1,...,h_d)=\Big(\begin{pmatrix} 1 & h_1\\ 0&1\end{pmatrix},...,\begin{pmatrix} 1& h_d\\0 &1\end{pmatrix}$\Big)}.$$

Since $\Gamma$ is a non-amenable subgroup of SL$_2(\mathbb Z)$, we get that $p_n:\Gamma\rightarrow X_n$ is onto for $n$ large enough. Since $\Lambda$ is abelian it does not have property $(\tau)$ with respect to $\{\ker(q_n)\}_{n=1}^{\infty}$. 
Thus,  Lemma \ref{tech2} provides an asymptotic homomorphism $\sigma_n:\Sigma\times\Lambda=(\Gamma*\mathbb Z)\times\Lambda\rightarrow\text{Sym}(X_n)$  such that
\begin{enumerate}
\item\label{uno} $\sigma_n(g,h)x=p_n(x)xq_n(h)^{-1}$, for all $g\in\Gamma,h\in\Lambda,x\in X_n$.
\item\label{dos} $\max\{\text{d}_{\text{H}}(\sigma_n(t,e)\circ\sigma_n(e,h),\sigma_n(e,h)\circ\sigma_n(t,e))\mid h\in\Lambda\}\geq\frac{1}{126}$, for infinitely many $n\geq 1$.
\end{enumerate}

Next, let $\widetilde X_n=\text{SL}_2(\mathbb Z/{r_{n,0}}\mathbb Z)\times X_n$ and define homomorphisms $\widetilde p_n:\Gamma\rightarrow\widetilde X_n$ and $\widetilde q_n:\Lambda\rightarrow\widetilde X_n$ by letting $\widetilde p_n(g)=(\pi_{r_{n,0}}(g),p_n(g))$ and $\widetilde q_n(h)=(e,q_n(h))$, for every $g\in\Gamma, h\in\Lambda$.
Further, we define  $\widetilde\sigma_n:\Sigma\times\Lambda\rightarrow\text{Sym}(\widetilde X_n)$ by letting for every $g\in\Sigma, h\in\Lambda, x\in \text{SL}_2(\mathbb Z/{r_{n,0}}\mathbb Z)$ and $y\in X_n$ \begin{equation}\label{defi}\text{$\widetilde\sigma_n(g,h)(x,y)=(\pi_{r_{n,0}}(g)x,\sigma_n(g,h)y)$.} \end{equation}
Then $(\widetilde\sigma_n)_{n\in\mathbb N}$ is an asymptotic homomorphism of $\Sigma\times\Lambda$ and conditions  \ref{uno} and \ref{dos} above rewrite as:
\begin{enumerate}[label=(\roman*)]
\item\label{i} $\widetilde\sigma_n(g,h)x=\widetilde p_n(g)x\widetilde q_n(h)^{-1}$, for all $g\in\Gamma,h\in\Lambda,x\in \widetilde X_n$.
\item\label{ii} $\max\{\text{d}_{\text{H}}(\widetilde\sigma_n(t,e)\circ\widetilde\sigma_n(e,h),\widetilde\sigma_n(e,h)\circ\widetilde\sigma_n(t,e))\mid h\in\Lambda\}\geq\frac{1}{126}$, for infinitely many $n\geq 1$.
\end{enumerate}
Since $\Gamma$ is a non-amenable subgroup of SL$_2(\mathbb Z)$, we get that $\widetilde p_n:\Gamma\rightarrow\widetilde X_n$ is onto for $n$ large enough. Moreover, a theorem of Bourgain and Varj\'{u} \cite[Theorem 1]{BV10} implies that $\Gamma$ has property $(\tau)$ with respect to $\{\ker(\widetilde p_n)\}_{n=1}^{\infty}$. By combining this fact with conditions (i) and (ii) above, we can apply Theorem \ref{conjugation} to conclude that there is no sequence of homomorphisms $\tau_n:\Sigma\times\Lambda\rightarrow\text{Sym}(Y_n)$, for some finite sets $Y_n\supset\widetilde X_n$, such that $\lim\limits_{n\rightarrow\infty}\text{d}_{\text{H}}(\widetilde\sigma_n(g),\tau_n(g)_{|\widetilde X_n})=0$, for every $g\in\Sigma\times\Lambda$.

Thus, in order to deduce that $\Sigma\times\Lambda$ is not weakly very flexibly P-stable, it suffices to argue that $(\widetilde\sigma_n)_{n\in\mathbb N}$ is a sofic approximation of $\Sigma\times\Lambda$. To see this, let $(g,h)\in(\Sigma\times\Lambda)\setminus\{(e,e)\}$.
 If $g\not=e$, then as $\lim\limits_{n\rightarrow\infty}r_{n,0}=+\infty$, we get that $\pi_{r_{n,0}}(g)\not=e$, for $n$ large enough. By using the definition \eqref{defi} of $\widetilde\sigma_n$, we get that $\text{d}_{\text{H}}(\widetilde\sigma_n(g,h),\text{Id}_{\widetilde X_n})=1$, for $n$ large enough. If $g=e$, then $h\not=e$ and since $\lim\limits_{n\rightarrow\infty}r_{n,i}=+\infty$, for all $1\leq i\leq d$, we get that $\widetilde q_n(h)\not=e$, for $n$ large enough. By using the definition \eqref{defi} of $\widetilde\sigma_n$, we get that $\text{d}_{\text{H}}(\widetilde\sigma_n(e,h),\text{Id}_{\widetilde X_n})=1$, for $n$ large enough. Since $\widetilde\sigma_n(e,e)=\text{Id}_{\widetilde X_n}$, for all $n\in\mathbb N$, this proves that $(\widetilde\sigma_n)_{n\in\mathbb N}$ is a sofic approximation of $\Sigma\times\Lambda$, finishing the proof of {\bf Case 1}.

{\bf Case 2.} $\Lambda=\mathbb F_k$, for some $k\geq 2$.

View $\Lambda$ as a subgroup of $\text{SL}_2(\mathbb Z)$ and let $\rho:\Lambda\rightarrow\text{SL}_2(\mathbb Z)$ be a homomorphism such that $\rho(\Lambda)\cong\mathbb Z$. For instance, if $a_1,...,a_k\in\Lambda$ are generators, we can let $\rho(a_1)=\begin{pmatrix}1&1\\0&1 \end{pmatrix}$ and $\rho(a_2)=...=\rho(a_k)=e$.

Fix $n\in\mathbb N$ and let $r_{n,0}, r_{n,1},r_{n,2}$ be $3$ distinct primes greater than $n$. Define $X_n=\prod_{i=1}^2\text{SL}_2(\mathbb Z/r_{n,i}\mathbb Z)$ and homomorphisms $p_n:\Gamma\rightarrow X_n, q_n:\Lambda\rightarrow X_n$ by letting for $g\in\Gamma$ and $h\in\Lambda$  $$\text{$p_n(g)=(\pi_{r_{n,1}}(g),\pi_{r_{n,2}}(g))$ and $q_n(h)=(\pi_{r_{n,1}}(\rho(h)),\pi_{r_{n,2}}(h)).$}$$

Since $\Gamma$ is a non-amenable subgroup of SL$_2(\mathbb Z)$, we get that $p_n:\Gamma\rightarrow X_n$ is onto for $n$ large enough. Since the image of $\rho$ is infinite abelian and $\lim\limits_{n\rightarrow\infty}r_{n,1}=+\infty$, $\Lambda$ does not have property $(\tau)$ with respect to $\{\ker(\pi_{r_{n,1}}\circ\rho)\}_{n=1}^{\infty}$. Since $\ker(q_n)\subset\ker(\pi_{r_{n,1}}\circ\rho)$, for every $n\in\mathbb N$, it follows that $\Lambda$ does not have property $(\tau)$ with respect $\{\ker(q_n)\}_{n=1}^{\infty}$.
Applying  Lemma \ref{tech2} provides an asymptotic homomorphism $\sigma_n:\Sigma\times\Lambda=(\Gamma*\mathbb Z)\times\Lambda\rightarrow\text{Sym}(X_n)$  which satisfies conditions \ref{uno} and \ref{dos} from above.

Next, let $\widetilde X_n=\text{SL}_2(\mathbb Z/{r_{n,0}}\mathbb Z)\times X_n$ and define the homomorphisms $\widetilde p_n:\Gamma\rightarrow\widetilde X_n$,  $\widetilde q_n:\Lambda\rightarrow\widetilde X_n$ and the asymptotic homomorphism  $\widetilde\sigma_n:\Sigma\times\Lambda\rightarrow\text{Sym}(\widetilde X_n)$
by the same formulae as in the proof of {\bf Case 1}. Then $(\widetilde\sigma_n)_{n\in\mathbb N}$ satisfies conditions \ref{i} and \ref{ii} from above. If $h\in\Lambda\setminus\{e\}$, then since $\lim\limits_{n\rightarrow\infty}r_{n,2}=+\infty$, we get that $\pi_{r_{n,2}}(h)\not=e$, for $n$ large enough. This implies that $\text{d}_{\text{H}}(\widetilde\sigma_n(e,h),\text{Id}_{\widetilde X_n})=1$, for $n$ large enough.  By repeating verbatim the rest of the argument from the proof of {\bf Case 1}, it follows that $\Sigma\times\Lambda$ is not weakly very flexibly P-stable. This finishes the proof of {\bf Case 2}. 

Finally,  the proof of Corollary \ref{B} shows that any group from parts (1)-(3) in its statement has a finite index subgroups which is isomorphic to either $\mathbb F_m\times\mathbb Z^d$ or to $\mathbb F_m\times\mathbb F_k$, for some $m,k\geq 2$ and $d\geq 1$. Thus, any group from Corollary \ref{B}, parts (1)-(3), is not weakly very flexibly P-stable.
\hfill$\blacksquare$

\end{document}